\documentclass[12pt,leqno]{article}
\usepackage{amsmath, amsthm, amsfonts, amssymb, color}
\usepackage{mathrsfs}
\theoremstyle{definition}

\newcommand{\scr}[1]{\mathscr #1}
\definecolor{wco}{rgb}{0.5,0.2,0.3}

\numberwithin{equation}{section} \theoremstyle{remark}

\newcommand{\ua}{\uparrow}

\title{
{\bf  Integrability Conditions for SDEs and Semi-Linear SPDEs
\footnote{Supported in part by NNSFC(11131003, 11431014), the 985 project.}}
}
\author{
{\bf Feng-Yu Wang  }\\
   \footnotesize{Laboratory of Mathematical and  Complex Systems, Beijing Normal University, Beijing 100875, China}\\
    \footnotesize{Department of Mathematics, Swansea University, Singleton Park, SA2 8PP, UK}\\
\footnotesize{Email: \tttext{wangfy@bnu.edu.cn}; \tttext{F.Y.Wang@swansea.ac.uk}} }

\begin{document}
\def\tttext#1{{\normalfont\ttfamily#1}} \def\div{{\rm div}}
\def\R{\mathbb R}  \def\ff{\frac} \def\ss{\sqrt} \def\B{\mathbf B}
\def\N{\mathbb N} \def\kk{\kappa} \def\m{{\bf m}}
\def\dd{\delta} \def\DD{\Delta} \def\vv{\varepsilon} \def\rr{\rho}
\def\<{\langle} \def\>{\rangle} \def\GG{\Gamma} \def\gam{\gamma}
  \def\nn{\nabla} \def\pp{\partial} \def\EE{\scr E}
\def\d{\text{\rm{d}}} \def\bb{\beta} \def\aa{\alpha} \def\D{\scr D}
  \def\si{\sigma} \def\ess{\text{\rm{ess}}}
\def\beg{\begin} \def\beq{\begin{equation}}  \def\F{\scr F}
\def\Ric{\text{\rm{Ric}}} \def\Hess{\text{\rm{Hess}}}
\def\e{\text{\rm{e}}} \def\ua{\underline a} \def\OO{\Omega}  \def\oo{\omega}
 \def\tt{\tilde} \def\Ric{\text{\rm{Ric}}}
\def\cut{\text{\rm{cut}}} \def\P{\mathbb P}
\def\C{\scr C}     \def\E{\mathbb E}
\def\Z{\mathbb Z} \def\II{\mathbb I}
  \def\Q{\mathbb Q}  \def\LL{\Lambda}\def\L{\scr L}
  \def\B{\scr B}    \def\ll{\lambda}
\def\vp{\varphi}\def\H{\mathbb H}\def\ee{\mathbf e}
\def\O{\scr O}

\maketitle
\begin{abstract} By using the  local  dimension-free Harnack inequality established on incomplete Riemannian manifolds, integrability conditions on the coefficients are presented for SDEs to imply the non-explosion of solutions  as well as the  existence, uniqueness and regularity estimates of invariant probability measures. These conditions include a class of drifts unbounded on compact domains such that the usual Lyapunov conditions can not be verified. The main results are extended to second order differential operators   on Hilbert spaces and  semi-linear SPDEs.
\end{abstract} \noindent

 \vskip 2cm

\section{Introduction}

In recent years, the existence and uniqueness   of strong solutions up to life time have been proved under  local integrability conditions for non-degenerate SDEs, see \cite{BC,KR,GM,Zhang,Zhang2} and references within.  See also \cite{CDR,DFRV,DR,W15,WZa,WZb} for extensions to degenerate SDEs and semi-linear SPDEs.

As a further development in this direction,   the present paper provides reasonable integrability conditions for the non-explosion of solutions, as well as the existence, uniqueness and regularity estimates of invariant probability measures. An essentially new point in the study is to make use of a local Harnack inequality in the spirit of \cite{W97}. With this inequality we are able to prove the non-explosion of a weak solution constructed from the Girsanov transform, see the proof of Lemma 3.1 below for details. Moreover, we use the hypercontractivity of the reference Markov semigroup to prove the boundedness of a Feyman-Kac semigroup induced by the singular SDE under study, which enables us to prove the existence of the invariant probability measure as well as a formula for the derivative of the density, see \eqref{Q2} and the proof of Lemma \ref{LNN} below for details.  To explain the motivation of the study more clearly, below we first recall some existing results in the literature,  then present a simple example to show how far can we go beyond.

Let $W_t$ be the $d$-dimensional Brownian motion on a complete filtration probability space $(\OO, \F, \{\F_t\}_{t\ge 0}, \P)$. Consider the following SDE (stochastic differential equation)   on $\R^d$:
\beq\label{SDE}\d X_t = b(X_t)\d t +\ss 2\,\si(X_t)\d W_t,\end{equation}  where $b: \R^d\to \R^d$ is measurable and
$\si\in W_{loc}^{1,1}(\R^d\to \R^d\otimes\R^d;\d x)$ such that $\si(x)$ is invertible for every $x\in\R^d.$
According to \cite[Theorem 1.1]{Zhang2} (see also  \cite{BC,KR,GM,Zhang} for earlier results), if  $|b|+\|\nn\si\|\in L^{p}_{loc}(\d x)$ for some $p>d$,
 then for any initial point $x$ the SDE \eqref{SDE} has a unique solution $(X_t^x)_{t\in [0, \zeta^x)}$ up to life time $\zeta^x$. We note that in \cite{Zhang2}   the global integrability and the uniform ellipticity conditions are assumed,   but  these conditions can be localized   since   for the existence and   uniqueness up to life time one only needs to consider   solutions  before   exiting bounded domains. On the other hand,   the ODE
 $$\d X_t= b(X_t)\d t$$
 does not have  pathwise uniqueness if $b$ is   merely H\"older continuous (for instance, $d=1$ and $b(x):=|x|^\aa$ for some $\aa\in (0,1)$).  So, the above result on SDE  indicates that the Brownian noise may $``$regularize" the drift to make an ill-posed equation   well-posed.

Next, sufficient integrability conditions for the non-explosion  have also been presented in \cite{Zhang}. For instance, if $\si$ is bounded and \beq\label{*01} |b|\le C+F\ \text{ for\ some\ constant}\  C>0  \text{\ and}\   F\in L^p(\d x)  \ \text{for\  some\ }  p>d,\end{equation} then the solution to \eqref{SDE}  is non-explosive. As the Lebesgue measure is infinite, this   condition is very restrictive. So, one of our aims is to replace it by   integrability conditions with respect to a probability measure, see Theorem \ref{T1.0} and Corollary \ref{C2.2} below.

We would like to indicate that when the invariant measure $\mu$ is given, there exist criteria on the conservativeness of non-symmetric Dirichlet forms,   which imply the non-explosion of solutions for $\mu$-a.e.   initial points, see \cite{TT} and references within. However, in our study the invariant probability measure is unknown, which is indeed the main object to characterize.
In general, to prove the existence of invariant probability measures   one uses Lyapunov (or drift) conditions. For instance,
if there exists  a positive   function $W_1\in C^2(\R^d)$ and a positive compact function $W_2$ such that
\beq\label{CP}LW_1:=\sum_{i,j=1}^\infty (\si\si^*)_{ij} \pp_i\pp_j W_1+ \sum_{i=1}^d b_i \pp_i W_1\le C-W_2 \end{equation} holds for some constant $C>0$,  then the associated diffusion semigroup has an invariant probability measure $\mu$ with $\mu(W_2)\le C$, see for instances \cite{H80, BKR, [2], BRW}.  Obviously, this condition is not available when $b$ is unbounded on compact sets. Our second purpose is to present a reasonable integrability condition for the existence and uniqueness of invariant probability measures, which applies to a class of SDEs with locally unbounded coefficients.  

Moreover, we also intend to investigate the  regularity properties of the invariant probability measure.   Recall that a probability measure $\mu$ on $\R^d$ is called an invariant probability measure of the generator $L$ (denoted by $L^*\mu=0$), if
\beq\label{2.2} \mu(Lf):= \int_{\R^d} Lf \d\mu=0,\ \ f\in C_0^\infty(\R^d).\end{equation} Obviously, an invariant probability measure $\mu$ of the  Markov semigroup $P_t$ associated to \eqref{SDE}  satisfies  $L^*\mu=0$.
In the past two decades, the existence, uniqueness and regularity estimates for invariant  probability measures of $L$ have been intensively investigated in both finite and infinite dimensional spaces, see the survey paper \cite{BKR} for concrete results and historical remarks. Here, we would like to recall a fundamental result on the regularity of the invariant probability measures.
Let $W_{loc}^{1,1}(\d x)$    be the class of  functions $f\in L_{loc}^1(\d x)$ such that
 $$\int_{\R^d}  f(x) (\div G)(x)\d x=-\int_{\R^d} \<G, F\>(x)\d x,\ \ G\in C_0^\infty(\R^d\to\R^d) $$ holds for some $F\in L^1_{loc}( \R^d\to \R^d;\d x)$, which is called   the weak gradient of $f$ and is denoted by  $F=\nn f$ as in the classical case.  For any $p\ge 1$, let
 $$W^{p,1}(\d x)=\big\{f\in W_{loc}^{1,1}(\d x):\ f,|\nn f|\in L^p(\d x)\big\}.$$
  Consider   the elliptic differential operator $L:=\DD+b\cdot \nn$ on $\R^d$ for some locally integrable $b:\R^d\to \R^d$.
  It has been shown in \cite{[5]}
that any invariant probability measure $\mu$ of $L$ with $\mu(|b|^2):=\int_{\R^d} |b|^2\d\mu<\infty$ has a density $\rr:=\ff{\d\mu}{\d x}$
such that $\sqrt \rr \in W^{2,1}(\d x)$. In addition,
\begin{equation}\label{1.2}
\int_{\mathbb{R}^d} \big|\nabla \sqrt \rr\big|^2\, \d x
\le \frac{1}{4} \int_{\mathbb{R}^d}|b|^2\, \d\mu.
\end{equation}

Since the invariant probability measure $\mu$ of $L$ is in general unknown, the integrability condition $\mu(|b|^2)<\infty$ is not explicit.
As   mentioned above that to ensure the existence of  $\mu$  one uses the Lyapunov condition \eqref{CP} for some positive function $W_1\in C^2(\R^d)$ and a compact function $W_2$, and   to verify $\mu(|b|^2)<\infty$ one would further need $|b|^2\le c +cW_2$   for some constant $c>0$. As we noticed  above that these   conditions do not apply if the coefficients   merely satisfy an  integrability condition  with respect to a   reference probability measure.

In conclusion,    we aim to search for   explicit integrability conditions on $b$ and $\si$ with respect to a  nice reference measure (for instance, the Gaussian measure) to  imply the non-explosion of solutions to the SDE \eqref{SDE}; the strong Feller property of the associated Markov semigroup;     the  existence, uniqueness and regularity estimates of the invariant probability measure. We also aim to extend the resulting assertions to the infinite-dimensional case.

The main results of this paper will be stated  in Section 2. Their proofs  are then presented in Sections 3-6 respectively. Finally, in Section 7 we present   a local  Harnack inequality  which plays a crucial role in the study.

To conclude this section, we present below a simple example to compare our results with existing ones introduced above.

\paragraph{Example 1.1.} Consider, for instance, the following SDE on $\R^d$:
$$\d X_t= \{Z(X_t)- \ll_0X_t\}\d t +\ss 2\,\d W_t,$$ where $\ll_0\in\R$ is a constant and $Z:\R^d\to \R^d$ is measurable.

(1) By Theorem \ref{T1.0} below for $\psi(x)=|x|$, if
\beq\label{XG1} \int_{\R^d} \e^{\vv |Z(x)|^2-\vv^{-1}|x|^2}\d x<\infty\ \text{for \ some\ }\vv\in (0,1),\end{equation}
then for any initial value the SDE has a unique strong solution which is non-explosive, and the associated Markov semigroup $P_t$ is strong Feller with a strictly positive density. Obviously, there are a lot   of maps $Z$ satisfying \eqref{XG1} but \eqref{*01} and the Lyapunov condition does not hold. For instance, it is the case when 
\beq\label{*002} Z(x):= x_0 \bigg\{\sum_{n=1}^\infty \log(1+|x-nx_0|^{-1})\bigg\}^{\theta}\end{equation}for some $x_0\in\R^d$ with $|x_0|=1$ and $\theta\in (0,\ff 1 2]$. 

(2) When $\ll_0>0$, we let $\mu_0(\d x)=  C\e^{-\ff{\ll_0}2|x|^2}\d x$ be a probability measure with normalization constant $C>0$. It is well known by Gross \cite{Gross} that  the log-Sobolev inequality in Assumption {\bf (H1)} holds for $\kk= \ff 2 {\ll_0}$ and $\bb=0$. By Theorem \ref{T1.1}, if
\beq\label{XG3} \int_{\R^d} \e^{\ll|Z(x)|^2-\ff{\ll_0 } 2|x|^2} \d x<\infty\ \text{for\ some\ }\ll>\ff 1{2\ll_0},\end{equation}
then $P_t$ has a unique invariant probability measure $\mu(\d x)=\rr(x)\d x$ such that
\beg{align*} &\mu_0\big(|\nn \ss\rr|^2\big)\le \ff{\ll_0}{4\ll\ll_0-2} \log\mu_0(\e^{\ll|Z|^2})<\infty,\\
&\mu_0\big(|\nn\log\rr\big|^2)\le \mu_0(|Z|^2)<\infty.\end{align*}
Obviously, for any $\theta\in (0,\ff 1 2)$, condition \eqref{XG3} holds for $Z$ defined by   \eqref{*002}, but  the Lyapunov condition \eqref{CP} is not available.

\section{Main results}

In the following four subsections, we introduce the main results in finite-dimensions and their infinite-dimensional extensions respectively.
To apply integrability conditions with respect to a reference measure $\mu_0$,  we regard the original SDE as a perturbation to the corresponding reference SDE whose semigroup is symmetric in $L^2(\mu_0)$.

\subsection{Non-explosion and strong Feller for SDEs}

Let $\si\in C^2(\R^d\to\R^d\otimes\R^d)$ with $\si(x)$ invertible for $x\in\R^d$ and denote $a= \si\si^*=(a_{ij})_{1\le i,j\le d}$.
For  $V\in C^2(\R^d)$,   define
\beq\label{L0}\beg{split} &Z_0= \sum_{i,j=1}^d \{\pp_j a_{ij}   - a_{ij} \pp_j V\}e_i,\\
& L_0={\rm tr}(a\nn^2)  +Z_0\cdot\nn =\sum_{i,j=1}^d a_{ij}\pp_i\pp_j +\sum_{i=1}^d
\<Z_0,e_i\>\pp_i,\end{split} \end{equation} where $\{e_i\}_{i=1}^d$ is the canonical orthonormal basis of $\R^d$, and $\pp_i$ is the directional derivative along $e_i$.

By  the integration by parts formula,   $L_0$ is symmetric in $L^2(\mu_0)$ for $\mu_0(\d x):= \e^{-V(x)}\d x:$
$$\mu_0(fL_0g)=-\mu_0(\<a\nn f,\nn g\>),\ \ f,g\in C_0^\infty(\R^d).$$ Then
$$\EE_0(f,g):=\mu_0(\<a\nn f,\nn g\>),\ \ f,g\in H^{2,1}_\si(\mu_0)$$ is a symmetric Dirichlet form generated by $L_0$, where $H^{2,1}_\si(\mu_0)$ is the closure of $C_0^\infty(\R^d)$ under the norm
  $$\|f\|_{H^{2,1}_\si(\mu_0)}:=\big\{\mu_0(|f|^2+|\si^*\nn f|^2)\big\}^{\ff 1 2}.$$  When $\si\equiv I$ (the identity matrix),  we simply denote
$H^{2,1}_\si(\mu_0)=H^{2,1}(\mu_0)$.

  Let $W_t$ be the $d$-dimensional Brownian motion as in Introduction. Consider  the reference SDE
\beq\label{E0} \d X_t= Z_0(X_t) \d t +\ss 2\,\si(X_t)\d W_t. \end{equation}
 Since   $\si$ and $Z_0$ are locally Lipschitz continuous, for  any initial point $x\in \R^d$ the SDE \eqref{E0} has a unique solution $X_t^x$ up to the explosion time $\zeta^x$. Let $P_t^0$ be the associated (sub-)Markov semigroup:
 $$ P_t^0f(x):= \E \big\{1_{\{\zeta^x>t\}} f(X_t^x)\big\},\ \ f\in \B_b(\R^d), t\ge 0, x\in \R^d.$$
When $\mu_0(\d x):= \e^{-V(x)}\d x$ is finite and   $1\in H^{2,1}_\si(\mu_0)$ with $\EE_0(1,1)=0$, we have   $P_t^01=1\  \mu_0$-a.e. Since
 $P_t^0 1$ is continuous  (indeed, differentiable)  for $t>0$, we have $P_t^01(x)=1$ for all $t\ge 0$ and $x\in\R^d.$
 Therefore, in this case the solution to \eqref{E0} is non-explosive for any initial points. By the symmetry of $P_t^0$ in $L^2(\mu_0)$,
 $\mu_0$ is $P_t^0$-invariant.

Now, for a measurable drift $Z:\R^d\to \R^d$, we consider the perturbed SDE
\beq\label{E1} \d X_t= \big\{Z+ Z_0\big\}(X_t)\d t +\ss 2\,\si(X_t)\d W_t. \end{equation}
By It\^o's formula, the generator of the solution is   $L:=L_0+Z\cdot\nn$.
According to \cite[Theorem 1.1]{Zhang2}, if $|Z|\in L^p_{loc}(\d x)$ for some $p>d$, then for any initial point $x\in \R^d$, the SDE \eqref{E1} has a  unique solution $X_t^x$ up to the life time
$\zeta^x$. We let $P_t$ be the associated (Dirichlet) semigroup:
$$P_t f(x)= \E\big[1_{\{t<\zeta^x\}} f(X_t^x)\big],\ \  x\in \R^d, t\ge 0, f\in \B_b(\R^d).$$
If $\P(\zeta^x=\infty)=1$ for all $x\in \R^d$, the solution is called non-explosive. In this case $P_t$ is a Markov semigroup.
More generally, for any non-empty open set $\O\subset \R^d$, let
$$T_\O^x= \zeta^x\land \inf\{t\in [0,\zeta^x): X_t^x \notin \O\},\ \ \inf\emptyset:=\infty.$$ Then the associated Dirichlet semigroup on $\O$ is given by
$$P_t^\O f(x)= \E\big[1_{\{t<T_\O^x\}} f(X_t^x)\big],\ \  x\in \O, t\ge 0, f\in \B_b(\O).$$

Let $\rr_\si$ be the  intrinsic metric induced by $\si$ as follows:
$$\rr_\si(x,y):=\sup\big\{|f(x)-f(y)|:\ f\in C^\infty(\R^d), |\si^*\nn f|\le 1\big\},\ \ x,y\in \R^d.$$
We have the following result.

\beg{thm}\label{T1.0} Let $\si\in C^2(\R^d\to\R^d\otimes\R^d)$ with $\si(x)$ invertible for $x\in\R^d$, and let $V\in C^2(\R^d)$ such that
\beq\label{PQ2}  \int_{\R^d} \Big( |\si^*\nn  \psi(x)|^2 +\e^{\vv |(\si^{-1}Z)(x)|^2}\Big)\e^{-V(x)-\vv^{-1}\rr_\si(0,x)^2  } \d x<\infty \end{equation}  holds for some constant $\vv\in (0,1)$ and a local Lipschitz continuous compact function $\psi$ on $\R^d$.  Then $\eqref{E1}$  has a unique  non-explosive solution for any initial points, and the associated Markov semigroup  $P_t$  is strong Feller with at most one invariant probability measure. Moreover, for any non-empty open set $\O\subset \R^d$ and $t>0$, $P_t^\O$ is strong Feller and has a strictly positive density $p_t^\O$   with respect to the Lebesgue measure on $\O$.   \end{thm}

\paragraph{Remark 2.1.} (1) Typical choices of $\psi$ include $|x|, \log (1 +|x|), \log\log (\e +|x|)...$    For instance, with   $\psi(x):= \log\log (\e+|x|)$  one may replace the term $|\si^*\nn  \psi(x)|^2$ in \eqref{PQ2} by $\ff{\|\si(x)\|^2}{(\e+|x|)^2\{\log (\e +|x|)\}^2}.$ So, if $V=0$ and $\int_{\R^d} \e^{-\ll \rr_\si(0,\cdot)^2}\d x<\infty$ for some $\ll>0$, the condition \eqref{PQ2} holds provided
$$ \log\ff{\|\si\|^2}{(\e+|\cdot|)^2\{\log(\e+|\cdot|)\}^2}+ |\si^{-1}Z|^2\le C (1+\rr_\si(0,\cdot))^2 + f$$   for some  constant $C>0$ and some function $f$ with $\e^{\vv f - \vv^{-1} \rr_\si(0,\cdot)^2}\in L^1(\d x)$ for some $\vv\in (0,1)$.

 (2) Let $\bar\si(r)=\sup_{|x|=r} \|\si(x)\|$ for $ r\ge 0$.   Then
$$\rr_\si(0,x)\ge U(x):=  \int_0^{|x|}\ff {\d r} {\bar\si(r)},\ \ x\in\R^d.$$ So, in \eqref{PQ2} we may replace $\rr_\si(0,\cdot)$ by the more explicit function $U$.

(3) The condition $\si\in C^2(\R^d\to\R^d\otimes\R^d)$ is stronger than $ \si \in  W^{p,1}_{loc}(\R^d\to\R^d\otimes\R^d;\d x)$ for some $p>1$ as required for the existence and uniqueness of solutions according to \cite[Theorem 1.1]{Zhang2}. This stronger condition is introduced because it together with the invertibility of $\si$ implies the local Harnack inequality  (see Theorem \ref{TA} below), which is a crucial tool in our study. If the local Harnack inequality could be established under weaker conditions, this condition would be weakened automatically. Indeed, under an additional assumption,  this condition will be replaced by $ \si \in  W^{p,1}_{loc}(\R^d\to\R^d\otimes\R^d;\d x)$ for some $p>1$,   see Theorem \ref{T1.1'} below for details.

\

Intuitively,  the non-explosion is a long distance property of the solution. So, it is natural for us to weaken the integrability condition \eqref{PQ2} by taking the integral outside a compact set.  But under this weaker condition we are not able to prove other properties included in Theorem \ref{T1.0}.

  \beg{cor} \label{C2.2} Let $\si\in C^2(\R^d\to\R^d\otimes\R^d)$ with $\si(x)$ invertible for $x\in\R^d$, and let $V\in C^2(\R^d)$ such that
\beq\label{PQ*}  \int_{D^c} \Big( |\si^*\nn  \psi(x)|^2 +\e^{\vv |(\si^{-1}Z)(x)|^2}\Big)\e^{-V(x)-\vv^{-1}\rr_\si(0,x)^2  } \d x<\infty \end{equation}  holds for some compact set $D\subset \R^d$, some  constant $\vv\in (0,1)$,  and some local Lipschitz continuous compact function $\psi$ on $\R^d$. If $Z\in L_{loc}^p(\d x)$ for some constant $p>d$, then the SDE $\eqref{E1}$ has a unique non-explosive solution for any initial points. \end{cor}

\subsection{Invariant probability measure for SDEs}

To investigate the invariant probability measures for the SDE \eqref{E1}, we  need the non-explosion of solutions such that
the standard tightness argument for the existence of invariant probability measure applies. To this end, we will apply Theorem \ref{T1.0} above, for which  we first assume that $\si$ is $C^2$-smooth (see {\bf (H1)} below)  then extend to less regular $\si$ by approximations (see {\bf (H$1'$)} below).

\ \newline {\bf Assumption (H1)}
\beg{enumerate} \item[(1)]  $\si\in C^2(\R^d\to\R^d\otimes\R^d)$ with $\si(x)$ invertible for $x\in\R^d$,    $V\in C^2(\R^d)$  such that $\mu_0(\d x):=\e^{-V(x)}\d x$ is  a probability measure satisfying
 \beq\label{SOB}  H^{2,1}_{\si}(\mu_0)=W^{2,1}_{\si}(\mu_0):=\big\{f\in W_{loc}^{1,1}(\d x):\ f,|\si^* \nn f|\in L^2(\mu_0)\big\}.\end{equation}   \item[(2)] The (defective) log-Sobolev inequality
\beq\label{LS} \mu_0(f^2\log f^2)\le \kk \mu_0(|\si^*\nn f|^2)+\bb,\ \ f\in C_0^\infty(\R^d), \mu_0(f^2)=1 \end{equation}  holds   for some   constants $\kk>0,\bb\ge 0.$
  \end{enumerate}

 Since $\mu_0(\d x):= \e^{-V(x)}\d x$ is finite, \eqref{SOB} implies   $1\in H^{2,1}_\si(\mu_0)$ with $\EE_0(1,1)=0$, so that the solution to \eqref{E0} is non-explosive as explained above. We note that \eqref{SOB} holds if  the metric $\rr_{\si}$ is complete.   Indeed,
in this case the function $\rr_{\si}(0,\cdot)$ is  compact with $|\si^*\nn \rr_\si(0,\cdot)|=1$, so that for any $f\in W^{2,1}_\si(\mu_0)$ we have
$f_n:= f\{1\land (n+1-\rr_{\si}(0,\cdot))^+\}\in H_{\si}^{2,1}(\mu_0)$ for $n\ge 1$, and it is easy to see that $f_n\to f$ in the norm $\|\cdot\|_{H^{2,1}_{\si}(\mu_0)}.$

There are   plentiful sufficient conditions for the log-Sobolev inequality \eqref{LS} to hold. For instance, if $\si\si^*\ge \aa I$ and $\Hess_V\ge K I$ for some constants $\aa,K>0$, then the Bakry-Emery criterion \cite{BE} implies \eqref{LS} for $\kk= \ff{2}{K\aa}$.  In the case that $K$ is not positive,   the log-Sobolev inequality holds for some
  constant $\kk>0$ if $\mu_0(\e^{\ll|\cdot|^2})<\infty $ for some $\vv>- \ff K 2,$ see \cite[Theorem 1.1]{W01}. See also \cite{CGW} for  a Lyapunov type sufficient condition of the log-Sobolev inequality.

\beg{thm} \label{T1.1}  Assume {\bf (H1)} and that
\beq\label{PP2} \mu_0(\e^{\ll |\si^{-1}Z|^2}):=\int_{\R^d} \e^{\ll |\si^{-1}Z|^2}\d\mu_0<\infty \end{equation} holds for some constant $\ll>\ff \kk{4}.$ Let $P_t$ be the semigroup associated to $\eqref{E1}$, and let  $L=L_0+Z\cdot \nn$ for $L_0$ in $\eqref{L0}$.  Then: \beg{enumerate}  \item[$(1)$]  $L$ has an invariant probability measure $\mu$, which is absolutely continuous with respect to  $\mu_0$ such that the density  function  $\rr:=\ff{\d\mu}{\d\mu_0}$
is strictly positive with   $\ss\rr,\,\log\rr \in H^{2,1}_\si(\mu_0)$ and
\beq\label{PD}\mu_0\big(|\si^*\nn\ss{\rr}|^2\big) \le \ff{1}{4   \ll- \kk} \big\{\log\mu_0(\e^{\ll|\si^{-1}Z|^2}) +\bb\big\}< \infty;\end{equation}
\beq\label{2.4} \mu_0(|\si^*\nn\log\rr|^2):=\lim_{\dd\downarrow 0} \int_{\R^d} \ff{|\si^*\nn \rr|^2}{(\rr+\dd)^2}\,\d\mu_0\le  \mu_0(|\si^{-1}Z|^2) <\infty. \end{equation}
 \item[$(2)$] The measure $\mu$ is the unique invariant probability measure of $L$ and   $P_t$  provided
\beq\label{MU0} \mu_0(\e^{\vv \|\si\|^2 })<\infty{\rm \ for \ some \ constant\  \vv>0}.\end{equation}     \end{enumerate} \end{thm}

\paragraph{Remark 2.2.} (1) Simply consider the case that  $\si=\si_0=I$.   If  $\Hess_V\ge K$  for some $K>0$, then   {\bf (H1)} holds  for $\kk= \ff 2K$ and $\bb=0$.   So, when  $\mu_0(\e^{\ll |Z|^2})<\infty$ holds for some  $\ll> \ff 1 {2K}$,   Theorem \ref{T1.1} implies that  $L$ and $P_t$ have a unique invariant probability measure $\mu$, which is absolutely continuous with respect to $\mu_0$, and the density function satisfies $\rr:=\ff{\d\mu}{\d\mu_0}$ satisfies $\ss\rr,\ \log\rr\in  H^{2,1}(\mu_0)$ with
$$\mu_0(|\nn\ss\rr|^2)\le \ff{K}{4K  \ll- 2} \log\mu_0(\e^{\ll |Z|^2})<\infty;\ \
 \mu_0(|\nn\log\rr|^2) \le   \mu_0(|Z|^2)  <\infty.$$

(2) Under {\bf (H1)}, if the super log-Sobolev inequality
$$\mu_0(f^2\log f^2)\le r \mu_0(|\si^*\nn f|^2)+\bb (r),\ \ r>0, f\in C_0^\infty(\R^d), \mu_0(f^2)=1 $$
holds for some $\bb: (0,\infty)\to (0,\infty)$,  then    Theorem \ref{T1.1} applies when  \eqref{PP2} holds    for some $\ll>0$, and in this case \eqref{PD} reduces to
$$\mu_0\big(|\si^*\nn\ss{\rr}|^2\big) \le \inf_{\ll>0, r\in (0, 4 \ll)} \ff{1}{4  \ll- r} \big\{\log\mu_0(\e^{\ll |\si^{-1}Z|^2}) +\bb(r)\big\}< \infty.$$
According to e.g. \cite[Theorems 2.1(1) and  2.3(2)]{RW03b} for $M=\R^d$,   the super log-Sobolev inequality holds  provided  $a\ge \aa I$ for some constant $\aa>0$ and $ \Hess_V$ is bounded below with $\mu(\e^{\ll|\cdot|^2})<\infty$ for any $\ll>0$.  In particular, it is the case when $a=I, V(x)= c_1+c_2|x|^p$ for some constants $c_1\in\R, c_2>0$ and $p>2$. See  \cite{DS,Gross,W05} and references within for more discussions on the super log-Sobolev inequality and the corresponding semigroup property.

(3) To illustrate  the sharpness of condition \eqref{PP2} for some $\ll>\ff \kk 4$, let us consider $\si=\si_0=I$ and  $V(x)=c+\ff 1 2|x|^2$ for some constant $c\in\R$,  so that  {\bf (H1)}  holds  for   $\kk= 2$ and $\bb=0$.
Let $Z(x)= r x =\ff r 2   \nn |\cdot|^2(x)$ for some constant $r\ge 0$. It is trivial that $L$ has  an invariant probability measure  if and only if $r<1,$ which is equivalent to $\mu_0(\e^{\ll|Z|^2})<\infty$ for some $\ll > \ff \kk{4}=\ff 1 2.$

\

 Now, we extend Theorem \ref{T1.1}  to less regular $\si$ by using the following assumption to replace
 {\bf (H1)}.

\ \newline
{\bf Assumption (H$1'$)}
\beg{enumerate} \item[(1)]  $\si\in W_{loc}^{p,1}(\R^d\to\R^d\otimes\R^d;\d x)$ for some $p>d$, $\si(x)$ is invertible for every $x\in\R^d$, and    $ a:=\si\si^*\ge \aa  I$ for some constants $\aa>0$.
\item[(2)] $V\in C^2(\R^d)$  such that $\mu_0(\d x):=\e^{-V(x)}\d x$ is  a probability measure satisfying \eqref{SOB}  and
 \beq\label{LS'} \mu_0(f^2\log f^2)\le \kk' \mu_0(|\nn f|^2) + \bb,\ \ f\in C_0^\infty(\R^d), \mu_0(f^2)=1 \end{equation}
   for some   constants $\kk'>0,\bb\ge 0.$
\item[(3)] There exists a constant $p>1$ such that $a_{ij}\in H^{2,1}(\mu_0)\cap L^{2p}(\mu_0)$ for any $1\le i,j\le d$ and  $|\nn V|\in L^{\ff{2p}{p-1}}(\mu_0).$
  \end{enumerate}

Let $L$ and $P_t$ be in Theorem \ref{T1.1} associated to the SDE \eqref{E1}.

\beg{thm}\label{T1.1'} Assume {\bf (H$1'$)} and let $\mu_0\big(\exp[\ll|Z|^2]\big)<\infty$ hold for some $\ll>\ff{\kk'}{4\aa^2}.$ Then
  $L$ and $P_t$ have a unique invariant probability measure $\mu(\d x):=\rr(x)\mu_0(\d x)$ for some   strictly positive function $\rr$  such that   $\ss\rr, \,\log\rr\in H^{2,1}(\mu_0)$ with
\beq\label{*D1}   \mu_0\big(|\nn\ss{\rr}|^2\big) \le \ff{1}{4\aa^2 \ll- \kk'} \big\{\log\mu_0\big(\e^{\ll|Z|^2}\big) +\bb\big\}< \infty;\end{equation}
\beq\label{*D2} \mu_0(|\nn\log\rr|^2):=\lim_{\dd\downarrow 0} \int_{\R^d} \ff{|\nn \rr|^2}{(\rr+\dd)^2}\,\d\mu_0\le  \ff 1 {\aa^2}\mu_0(|Z|^2) <\infty. \end{equation}
 \end{thm}

\subsection{Elliptic  differential operators on Hilbert spaces }

 We first consider the invariant probability measure of second order differential operators on a separable Hilbert space, then apply  to semi-linear SPDEs.
We will take a Gaussian measure   as the reference measure.

Let $(\H, \<\cdot,\cdot\>, |\cdot|)$ be a separable Hilbert space, let $(A,\D(A))$ be a positive definite self-adjoint operator on $\H$ having discrete spectrum with all eigenvalues $(0<)\ll_1\le \ll_2\le \cdots$ counting multiplicities  such that
\beq\label{EG} \sum_{i=1}^\infty \ll_i^{-1}<\infty. \end{equation}
 Let $\{e_i\}_{i\ge 1}$ be the corresponding eigenbasis of $A$.
Let $\mu_0$ be the Gaussian measure on $\H$ with covariance operator $A^{-1}$.
In coordinates with respect to the basis $\{e_i\}_{i\ge 1}$, we have
\beq\label{MU} \mu_0(\d x)= \prod_{i=1}^\infty \Big(\ff {\ss{\ll_i}}{\ss{2\pi }} \e^{-\ff {\ll_i x_i^2}{2}}\,\d x_i\Big),\ \ \ x_i:=\<x,e_i\>, i\ge 1.\end{equation}
For any $n\ge 1$, let $\H_n={\rm span}\{e_i: 1\le i\le n\}$  and define the probability measure
$$\mu_0^{(n)}(\d x)= \prod_{i=1}^n \Big(\ff {\ss{\ll_i}}{\ss{2\pi }} \e^{-\ff {\ll_i x_i^2}{2}}\,\d x_i\Big) \ \text{on}\   \H_n.$$  We have $\mu_0^{(n)} =\mu_0\circ\pi_n^{-1}$ for the orthogonal
projection $\pi_n: \H\to\H_n$.

Let   $(\L(\H),\|\cdot\|)$ be the space of  bounded linear operators on $\H$ with operator norm $\|\cdot\|,$ and let $\L_s(\H)$ be the class of all symmetric elements in $\L(\H)$.  For any $a\in \L_s(\H)$ let  $a_{ij}=\<a e_i,e_j\>$ for $i,j\ge 1$. We make the following assumption.

\ \newline {\bf Assumption (H2)}
\beg{enumerate}\item[(1)] $a_{ij}\in C^2(\H)$ for $i,j\ge 1$, and  $a\ge \aa I$ for some  constant $\aa>0$.
 \item[(2)] For   $n\ge 1$ and  $\si_n:=\ss{(a_{ij})_{1\le i,j\le n}},$\ $H_{\si_n}^{2,1}(\mu_0^{(n)})= W^{2,1}_{\si_n}(\mu_0^{(n)})$ holds.
 \item[(3)] For any $i,j\ge 1$, there exists $\vv_{ij}\in (0,1)$ such that
\beq\label{SII} \sup_{n\ge 1} \int_{\R^n}\exp\big[\vv_{ij}   |a_{ij}|^{1+\vv_{ij}}\big] \d\mu_0^{(n)} <\infty.\end{equation} \end{enumerate}
We note that
$$\int_{\R^n}\exp\big[\vv_{ij}   |a_{ij}|^{1+\vv_{ij}}\big] \d\mu_0^{(n)}=  \int_{\H}\exp\big[\vv_{ij}   |a_{ij} \circ\pi_n|^{1+\vv_{ij}}\big] \d\mu_0.$$

As mentioned above that $H_{\si_n}^{2,1}(\mu_0^{(n)})= W^{2,1}_{\si_n}(\mu_0^{(n)})$ is implied by the completeness of the metric on $\R^n$ induced by $\si_n,$ and the later holds   if  for any $i,j\ge 1$ there exists $\vv_{ij}>0$ such that
 \beq\label{SII0} |a_{ij}(x)|\le \ff 1 {\vv_{ij}} (1+|x|)^2,\ \ x\in\H.\end{equation}
  The condition  \eqref{SII} will be used for finite-dimensional approximations in the end of the proof of Theorem \ref{T4.1}(1) below. According to \eqref{EG} and the definition of $\mu_0^{(n)}$,  the conditions \eqref{SII} and \eqref{SII0}   hold  provided for any $i,j\ge 1$ there exists a constant $\vv_{ij}'\in (0,1)$ such that $|a_{ij}(x)|\le \ff 1 {\vv_{ij}'} (1+ |x|)^{\ff 2 {1+\vv_{ij}'}}.$

Let  $\pp_i$ be the directional derivative along $e_i, i\ge 1.$ For a   measurable drift $Z: \H\to \H$, consider   the operators
\beq\label{LLL}L:= L_0+Z\cdot\nn,\ \ \ L_0 := \sum_{i,j=1}^\infty \Big(a_{ij} \pp_i\pp_j     + \big\{\pp_j a_{ij}  - a_{ij} \ll_j \big\} \pp_i\Big), \end{equation}which are well   defined on
the class of smooth cylindrical functions with compact support:
 $$  \F C_0^\infty:=\big\{\H\ni x\mapsto f(\<x,e_1\>,\cdots, \<x,e_n\>): n\ge 1, f\in C_0^\infty(\R^n)\big\}.$$ It is easy to see that $L_0$ is symmetric in $L^2(\mu_0)$:
 $$\mu_0(fL_0g) =-\mu_0(\<a\nn f,\nn g\>),\ \ f,g\in \F C_0^\infty.$$
Let $H^{2,1}(\mu_0)$ be the completion of $\F C_0^\infty$ with respect to the inner product
$$\<f,g\>_{H^{2,1}(\mu_0)} := \mu_0(fg)+ \mu_0(\<\nn f,\nn g\>).$$
A probability measure $\mu$ on $\H$ is called an invariant probability measure of $L$  (denoted by $L^*\mu=0$), if for any $f\in \F C_0^\infty$ we have  $Lf \in L^1(\mu)$ and  $\mu(Lf)=0$.

\beg{thm}\label{T4.1}  Assume $\eqref{EG}$ and {\bf (H2)}.\beg{enumerate} \item[$(1)$] If   $\mu_0(\e^{\ll |Z|^2})<\infty$ for some $\ll>\ff 1 {2\ll_1\aa^2}$, then
  $L$ has an invariant  probability measure  $\mu$, which is absolutely continuous with respect to $\mu_0,$ and the density function $\rr:=\ff{\d\mu}{\d\mu_0}$ satisfies     $\ss\rr,\ \log\rr\in H^{2,1}(\mu_0)$ with
\beq\label{PD'}\mu_0\big( \big|\nn \ss\rr\big|^2\big) \le \ff{\ll_1}{4\aa^2\ll_1\ll- 2} \log\mu_0(\e^{\ll|Z|^2}) <\infty \end{equation}
  and
\beq\label{2.4'}\mu_0\big(\big|\nn \log {\rr}\big|^2\big):= \lim_{\dd\downarrow 0}  \mu_0\big( \big|\nn \log {\rr+\dd}\big|^2\big)\le   \ff{\mu_0(|Z|^2)}{\aa^2}<\infty.\end{equation}
\item[$(2)$] If moreover  $\|a\|_\infty<\infty$, then $(L, \F C_0^\infty)$ is closable in $L^1(\mu)$ and the closure generates a Markov $C_0$-semigroup $T_t$ on $L^1(\mu)$ with $\mu$ as an invariant probability. Moreover, there exists a standard Markov process $\{\bar \P_x\}_{x\in\H\cup\{\pp\}}$ on $\H\cup\{\pp\}$ which is continuous and  non-explosive   for $\EE^\mu$-q.e. $x$, such that the associated Markov semigroup $\bar P_t$ is a $\mu$-version of $T_t$; that is, $\bar P_t f=T_tf\ \mu$-a.e. for all $t\ge 0$ and $f\in \B_b(\H).$  \end{enumerate}
\end{thm}

For readers' convenience, we would like to recall here the notion of standard Markov process involved in Theorem \ref{T4.1}(2). Let $\pp$ be an extra point and extend the topology of $\H$ to $\H\cup\{\pp\}$ by letting the set $\{\pp\}$  open.
A family of probability measures  $\{\P_x\}_{x\in\H\cup\{\pp\}}$ on $$\OO:=\big\{\oo: [0,\infty)\to \H\cup\{\pp\}: \text{if}\ \oo_t=\pp\ \text{then}\ \oo_s=\pp\ {\rm for}\  s\ge t\big\}$$ equipped with the   $\si$-field $\F:=\si(\oo_t:t\ge 0)$   is called a standard  Markov process, if $\P_x(\oo_0=x)=1$ and the distribution $P_t(x,\d y)$ of $\OO\ni\oo\mapsto \oo_t$ under $\P_x$ gives rise to a  Markov transition kernel on $\H \cup\{\pp\}$. When the process is non-explosive, i.e. $$\P_x\big(\inf\{t\ge 0: \oo_t=\pp\}=\infty\big)=1,\ \ x\in\H,$$  the sub-family $\{\P_x\}_{x\in \H}$ is a standard Markov process on $\H$. In this case,
the process (or the associated Markov semigroup  $P_t$) is called Feller if $P_tC_b(\H)\subset C_b(\H)$ for all $t\ge 0$, and is called strong Feller if $P_t\B_b(\H)\subset C_b(\H)$ for all $t>0.$ If moreover $\P_x(C([0,\infty)\to \H))=1$ holds for all $x\in\H$, then the process is    continuous.

\

Next, we  extend Theorem \ref{T1.1'} to the infinite-dimensional case, for which we need the following assumption.

\ \newline{\bf Assumption (H$2'$)}
\beg{enumerate}\item[(1)]   $a\ge \aa I$ for some  constant $\aa>0$, and  for every $n\ge 1$ there exists a constant $p>n$ such that $a_n\in W_{loc}^{p,1}(\R^n\to\R^n\otimes\R^n; \d x).$    \item[(2)]  For any  $i,j\in\mathbb N$ there exists $\vv_{ij}\in (0,1)$ such that \eqref{SII}  and $\mu_0^{(n)}(|\nn a_{ij}\circ\pi_n|^{2+\vv_{ij}})<\infty$ hold for any $n\ge 1$. \end{enumerate}

\beg{thm}\label{T4.1'} Under $\eqref{EG}$ and {\bf (H$2'$)},    assertions $(1)$ and $(2)$ in Theorem $\ref{T4.1}$ hold. \end{thm}

\subsection{Semi-linear SPDEs}

We intend to investigate the existence, uniqueness and non-explosion of the SPDE corresponding to $L$ in \eqref{LLL}, and  to show that   the probability measure in Theorem \ref{T4.1} is the unique invariant probability measure of the associated  Markov semigroup. For technical reasons, we only consider the case that $a=I$,
 for which the corresponding SPDE reduces to the standard semi-linear SPDE
\beq\label{E2} \d X_t= \big\{Z(X_t)-AX_t  \big\}\d t +\ss 2\,\d W_t,\end{equation}
where $Z:\H\to \H$ is measurable,
$W_t$ is the cylindrical Brownian motion, i.e.
$$W_t= \sum_{i=1}^\infty \bb_t^i e_i,\ \ t\ge 0$$ for a sequence of independent one-dimensional Brownian motions $\{\bb_t^i\}_{i\ge 1}.$ An adapted continuous process $X_t$ on $\H$ is called a mild solution to \eqref{E2}, if
$$X_t= \e^{-A t} X_0+\int_0^t \e^{-A(t-s)}Z(X_s)\d s + \int_0^t \e^{-(t-s)A}\d W_s,\ \ t\ge 0.$$
We assume
\beg{enumerate} \item[{\bf (H3)}]  $\sum_{i=1}^\infty  \ff 1 {\ll_i^{\theta}} <\infty$ for some $\theta\in (0,1)$, and $\mu_0(\e^{\ll |Z|^2})<\infty$ for some constant $\ll>0$. \end{enumerate}

According to the recent paper \cite{DFRV}, {\bf (H3)} implies the existence and pathwise uniqueness of mild solutions to \eqref{E2}  for $\mu_0$-a.e. starting points.
Below we intend to prove the weak uniqueness of \eqref{E2} for any initial points.   A standard continuous Markov process on $\H$ is called a weak solution to \eqref{E2}, if it solves the martingale problem for $(L, \F C_0^\infty)$. In this case one may construct a cylindrical Brownian motion $W_t$ on the probability space $(C([0,\infty)\to\H);\F, \P_x)$, where $\F:=\si(\{\oo\mapsto\oo_t: t\ge 0\})$,  such that the coordinate process $X_t(\oo):=\oo_t$ is a mild solution to \eqref{E2} with $X_0=x$. See e.g.  \cite[Proposition IV.2.1]{IW} for the explanation in the finite-dimensional case, which works also in the present case as the cylindrical Brownian motion is determined by its finite-dimensional projections.

\beg{thm}\label{T4.2} Assume that {\bf (H3)} holds.     \beg{enumerate} \item[$(1)$] There exists a standard continuous Markov process $\{\P_x\}_{x\in \H}$ solving $\eqref{E2}$ weakly for every initial point, and    the associated Markov semigroup $P_t$ is strong Feller having a strictly positive density with respect to $\mu_0$.
\item[$(2)$] If $Z$ is bounded on bounded sets, then there exists a unique standard Markov process   solving $\eqref{E2}$ weakly for every initial point such that the associated Markov semigroup is Feller.
\item[$(3)$] If $Z$ is bounded on bounded sets and $\mu_0(\e^{\ll|Z|^2})<\infty$ holds for some  $\ll>\ff 1 {2\ll_1},$ then $P_t$ has a unique invariant probability measure $\mu$, which is absolutely continuous with respect to $\mu_0$ and  the  density function $\rr:=\ff{\d\mu}{\d\mu_0}$ is strictly positive with $\ss\rr, \log\rr\in H^{2,1}(\mu_0)$  such that estimates $\eqref{PD'}$ and $\eqref{2.4'}$ hold for $\aa=1$. \end{enumerate}  \end{thm}

\paragraph{Remark 2.3.}   Unlike in the finite-dimensional case where   $Z\in L^p_{loc}(\d x)$  for some $p>d$ implies the pathwise uniqueness of the solution for any initial points,  in the infinite-dimensional case  this is unknown without any continuity conditions on $Z$. It is shown in \cite{W15} (also for the multiplicative noise case) that if $Z$ is Dini continuous then the pathwise uniqueness holds for any initial points.

\section{Proof of Theorem \ref{T1.0}}
The main idea is to show that the solution to the reference SDE \eqref{E0} is a weak solution to \eqref{E1} under a weighted probability, so that the non-explosion of \eqref{E0} implies that of \eqref{E1}. To this end, we will apply the local Harnack inequality \eqref{HI} below to verify the Novikov condition for the Girsanov transform. To realize the idea, we
  first consider the case that
\beq\label{PQ2'}   \int_{\R^d}   \e^{\vv |(\si^{-1}Z)(x)|^2 -V(x) } \d x<\infty  \end{equation} holds for some $\vv>0$,
then reduce back to the original condition \eqref{PQ2}.

\beg{lem}\label{P1} Assume that   $\eqref{PQ2'}$ holds for some constant $\vv>0$ and $1\in H_\si^{2,1}(\mu_0)$ with $\EE_0(1,1)=0$, then all assertions in Theorem $\ref{T1.0}$ hold.  \end{lem}

\beg{proof} Obviously, \eqref{PQ2'} implies that $\mu_0(\d x):=\e^{-V(x)}\d x$ is a finite measure. Since the coefficients in \eqref{E0} is locally Lipschitz continuous, it is classical that the SDE has a unique solution up to the explosion time. Since $1\in H_\si^{2,1}(\mu_0)$ with $\EE_0(1,1)=0$, as explained after \eqref{E0} that the solution to
\eqref{E0} is non-explosive and $\mu_0$ is $P_t^0$-invariant.
Moreover, since the drift in \eqref{E1} is locally bounded, according to \cite{Zhang}, this SDE has a unique solution for any initial points.
So, it remains to show that the solution is non-explosive, and  the associated Markov  semigroup $P_t$  is strong Feller with at most one invariant probability measure.

 A crucial tool in the proof is the following local Harnack inequality. Consider $\R^d$ with the $C^2$-Riemannian metric
  $$\<u,v\>_\si:= \<\si\si^* u,v\>,\ \ u,v\in\R^d,$$ and let $\DD_\si, \nn_\si$ be the corresponding Laplace-Beltrami operator and the gradient operator. Then $L_0$ can be rewritten as
  $$L_0=\DD_\si +\nn_\si \bar V$$ for some $\bar V\in C^2(\R^d).$ Since the intrinsic distance $\rr_\si$   is locally equivalent to the Euclidean distance, according to Theorem \ref{TA} below, for any $p>1$ there exists  positive  $\Phi_p\in C(\R^d)$  such that
\beq\label{HI}  (P_t^0 f)^p(x) \le (P_t^0 f^p (y))\exp\bigg[ \Phi_p(x) \Big(1 + \ff{|x-y|^2}{1\land t}\Big)\bigg]
 \ \   x,y\in \R^d, \ |x-y|\le \ff 1 {\Phi_p(x)} \end{equation}
holds for all $t>0$ and $f\in \B_b^+(\R^d):=\{f\in\B_b(\R^d):\ f\ge 0\}.$

\

(a) {\bf Non-explosion}. It suffices to find out a constant $t_0>0$ such that for any initial points, the solution to \eqref{E1} is non-explosive before time $t_0$.  To this end, we construct a weak solution by using   the reference SDE \eqref{E0}.    We intend   to find out $t_0>0$ such that for any initial point $x$, the solution to \eqref{E0} for $X_0=x$ is a weak solution to \eqref{E1} for $t\in [0,t_0].$ So, by the weak uniqueness of   \eqref{E1}, which follows from the strong uniqueness, we conclude that the strong solution to \eqref{E1} is non-explosive before $t_0.$  To this end,  we    verify the Novikov condition
\beq\label{NV} \E \exp\bigg[ \ff 1 4 \int_0^{t_0}|(\si^{-1}  Z)(X_s)|^2\d s\bigg]<\infty,\ \ X_0=x\in\R^d,\end{equation}
so that $\Q:= \exp[\ff 1 {\ss 2} \int_0^{t_0}\<(\si^{-1}  Z)(X_s),\d W_s\>- \ff 1 4 \int_0^{t_0}|(\si^{-1}  Z)(X_s)|^2]\P$ is a probability measure. In this case,  by the Girsanov theorem,
$$\tt W_t:= W_t-\ff 1 {\ss 2}\int_0^t (\si^{-1}  Z)(X_s)\, \d s,\ \ t\in [0,t_0]$$ is a Brownian motion under $\Q$. Thus, rewriting  \eqref{E0} as
 $$\d X_t = (  Z+  Z_0)(X_t) +\ss 2\, \si(X_t)\d \tt W_t,\ \ t\in [0,t_0],$$ we see that $(X_t,\tt W_t)_{t\in [0,t_0]}$ is a weak solution to \eqref{E1}
  under the probability measure $\Q$.

  To prove \eqref{NV}, we use the Harnack inequality \eqref{HI} for $p=d+1$ to derive
\beg{align*} &\big\{\E \e^{\ll (|(\si^{-1}  Z)(X_t)|^2\land N)}\big\}^{d+1}= \big( P_t^0 \e^{\ll (|\si^{-1}  Z|^2\land N)} (x)\big)^{d+1}\\
  &\le   P_t^0 \e^{(d+1)\ll (|\si^{-1}  Z|^2\land N)}(y) \e^{\Phi_{d+1}(x) (1 + |x-y|^2/t)},\ \ t\in (0,1], N>0, |y-x|\le \ff 1 {\Phi_{d+1}(x)}.\end{align*}  Since $\mu_0$ is $P_t^0$-invariant,   for $B_{x,t}  := \big\{y: |y-x|\le \ff 1 {\Phi_{d+1}(x)}\land \ss t\big\}$ this implies
 \beg{align*} & \big\{\E \exp\big[\ll  (|(\si^{-1}  Z)(X_t)|^2\land N)\big]\big\}^{d+1}   \mu_0 ( B_{x,t})
  \e^{-2\Phi_{d+1}(x)}\\
 &\le \int_{ B_{x,t}}\big(  P_t^0 \e^{ \ll   (|\si^{-1}  Z |^2\land N)}\big)^{d+1}(x)\, \exp\Big[-\Phi_{d+1}(x)\Big(1+\ff{|x-y|^2}{t}\Big)\Big]  \mu_0(\d y) \\
 &\le \int_{ B_{x,t}}  P_t^0 \e^{(d+1)\ll (|\si^{-1}  Z |^2\land N)}(y) \mu_0(\d y)\le  \mu_0( \e^{\vv |\si^{-1}  Z|^2})<\infty,\ \ t\in (0,1],
 \ll\in \Big(0,\ff\vv{d+1}\Big].\end{align*} Since $ \mu_0$ has strictly positive and continuous density $\e^{-V}$ with respect to $\d x$,  there exists  $G\in C(\R^d\to (0,\infty))$ such that $ \mu_0( B_{x,t})\ge G(x) t^{\ff d 2}$   for $t\in (0,1]$ and $x\in\R^d$.
 By taking $\ll=\vv/(d+1)$ and letting $N\to\infty$ in the above display, we arrive at
$$\E \e^{\ff\vv{d+1}  |(\si^{-1}  Z)(X_t)|^2} \le \ff{H(x)}{\ss t}<\infty,\ \ t\in (0,1], x\in \R^d$$ for some positive  $H\in C(\R^d).$
 Therefore, by Jensen's inequality,   we have
 \beq\label{EST2} \beg{split} &\E \exp\bigg[\gam \int_0^{r} |(\si^{-1}  Z)(X_s)|^2\d s\bigg] \le \ff 1 {r}\int_0^{r}
 \E\e^{\gam r |(\si^{-1}  Z)(X_s)|^2}\d s   \\
 &\le \ff 1 {r} \int_0^{r} \ff{H(x)}{\ss t}\d t= \ff{2H(x)}{\ss{r}}<\infty,\ \ x\in \R^d,   r\in (0,1], \gam\in   \Big(0,  \ff {\vv}{(d+1)r}\Big].\end{split}\end{equation}
 This implies \eqref{NV} by taking $\gam=\ff 14$ and  $t_0=r= 1\land \ff{4\vv}{d+1}.$

\

(b) {\bf Strong Feller of $P_t$ and uniqueness of invariant probability measure}. According to \cite[Theorem 4.1]{BKR0}, the Markov semigroup  $P_t^0$ is strong Feller. For any $x\in\R^d$, we let $X_s^x$ solve \eqref{E0} with initial point $x$ and define
$$R_r^x=\exp\bigg[\ff 1 {\ss 2}\int_0^r\<(\si^{-1}Z)(X_s^x),\d W_s\> -\ff 1 4 \int_0^r |(\si^{-1}Z)(X_s^x)|^2\d s\bigg],\ \ r\in [0,t_0].$$
By \eqref{NV} and the Girsanov theorem, we have
$$P_t f(x)= \E \big[f(X_t^x) R_t^x\big],\ \ t\in [0,t_0], f\in B_b(\R^d), x\in \R^d.$$ Then for any $t>0, x\in \R^d$ and $f\in \B_b(\R^d)$,
the semigroup property of $P_s$ and the strong Feller property of $P_s^0$ imply
\beg{align*}  & \limsup_{y\to x} |P_t f(y)-P_t f(x)| = \limsup_{y\to x} |P_r(P_{t -r}f)(y)-P_r(P_{t-r} f)(x)|\\
&= \limsup_{y\to x} \big|\E[R_r^y P_{t -r}f)(X_{r}^y)-R_r^x (P_{t-r} f)(X_r^x)]\big|\\
&\le \limsup_{y\to x} \Big\{ \big|P_r^0(P_{t -r}f)(y)- P_r^0(P_{t -r}f)(x)\big|+   \E\big(|R_r^y-1|+  |R_r^x -1|\big)\Big\}\\
&\le \sup_{y: |y-x|\le 1}  \E\big(|R_r^y-1|+  |R_r^x -1|\big),\ \ r\in (0, t).\end{align*}
Noting that $\E |R_r^y-1|^2= \E (R_r^y)^2 -1$ for small $r>0$, then the strong Feller property follows provided
\beq\label{STF} \limsup_{r\to 0} \sup_{y: |y-x|\le 1}\E(R_r^y)^2 \le 1.\end{equation} To prove this,
we let  $M_r = \ff 1 {\ss 2}\int_0^r\<(\si^{-1}Z)(X_s^x),\d W_s\>.$ Then for small $r>0$
$$ \E (R_r^y)^2 =\E\e^{2M_r -\<M\>_r} \le \big(\E\e^{4 M_r - 8\<M\>_r}\big)^{\ff 1 2}\big(\E\e^{6\<M\>_r}\big)^{\ff 1 2}
= \big(\E \e^{3\int_0^r |(\si^{-1}Z)(X_s^x)|^2\d s}\big)^{\ff 1 2}.$$  So, applying \eqref{EST2} with   $\gam=\ff\vv{(d+1)r}$ for small $r>0$, and using Jensen's inequality, we obtain
 \beg{align*} &\limsup_{r\to 0} \sup_{y: |y-x|\le 1}\big\{ \E (R_r^y)^2\big\}^{2} \le \limsup_{r\to 0}\sup_{y: |y-x|\le 1}\E \e^{3\int_0^r |(\si^{-1}Z)(X_s^x)|^2\d s}\\
 &\le\limsup_{r\to 0}\sup_{y: |y-x|\le 1}\big(\E \e^{\gam\int_0^r |(\si^{-1}Z)(X_s^x)|^2\d s}\big)^{\ff 3\gam}
  \le \limsup_{r\to 0}\sup_{y: |y-x|\le 1}\Big(\ff{2H(y)}{\ss r}\Big)^{\ff{3(d+1)r}\vv}=1.\end{align*} This implies \eqref{STF}.

Next, as already mentioned above,    every invariant probability measure of $P_t$  has strictly positive density with respect to  the Lebesgue measure, so that any two invariant probability measures are equivalent each other.   Therefore, the invariant probability measure has to be unique, since it is well known that any two different extremal   invariant probability measures of a strong Feller Markov operator are singular each other.

(c) {\bf The assertion for $P_t^\O$}. Due to the semigroup property ensured by the pathwise uniqueness, it suffices to prove for small enough $t>0$. Let $T_\O^x$ be the hitting time of $X_t^x$ to the boundary of $\O$. By the Girsanov theorem we have
\beq\label{SEM} P_t^\O f(x)= \E\big[1_{\{T_\O^x>t\}} f(X_t^x) R_t^x\big],\ \ f\in \B_b(\O), x\in \O.\end{equation}
Let $P_t^{\O,0}f(x)= \E\big[1_{\{T_\O^x>t\}} f(X_t^x)\big]$ be the Dirichlet semigroup associated to \eqref{E1}. Since $\si$ is invertible, by the $C^2$-regularity of $\si$ and $V$  we see that $P_t^{\O,0}$ is strong Feller having strictly positive density with respect to the Lebesgue measure
(see \cite{ATW13} for gradient estimates and log-Harnack inequalities of $P_t^{\O,0}$). Then the strong Feller property can be proved as in (b) using $P_t^{\O,0}$ in place of $P_t^0$.

Next, by \eqref{EST2} we have  $\E \{(R_t^x)^{-1}\}<\infty$ for small $t>0$. Then for any measurable set $A$ such that $P_t^\O1_A(x)=0$,    \eqref{SEM} implies
$$ \{P_t^{\O,0}1_A(x)\}^2= \big\{\E\big[1_{\{T_\O^x>t\}}  1_A (X_t^x)\big]\big\}^2\le \{P_t^{\O}1_A(x)\} \E\{(R_t^x)^{-1}\}=0.$$
Thus, the measure $P_t^{\O,0}1_{\d z}(x)$ is absolutely continuous to the measure $P_t^{\O}1_{\d z}(x).$ Since $P_t^{\O,0}$ has a strictly positive density, so does $P_t^{\O}.$  \end{proof}

\beg{proof}[Proof of Theorem \ref{T1.0}] Since $|\si^*\nn\rr_\si(0,\cdot)|=1$, for any $\dd>0$ the function $\rr_\si(0,\cdot)$ can be uniformly approximated by smooth ones $f_n$ with $|\si^*\nn f_n| \le 1+\dd$. In particular, we may take $\tt\rr\in C^2(\R^d)$ such that $|\rr_\si(0,\cdot)-\tt\rr|\le 1$ and $|\si^* \nn \tt \rr|^2\le 2$, so that \eqref{PQ2} holds for some $\vv\in(0,1)$ if and only if
\beq\label{PQ2''} \int_{\R^d} \Big( |\si^*\nn  \psi(x)|^2 +\e^{\vv |(\si^{-1}Z)(x)|^2}\Big)\e^{-V(x)-\vv^{-1}\tt\rr(x)^2  } \d x<\infty\end{equation} holds for some $\vv\in (0,1).$

 To apply Lemma \ref{P1}, we take
$$\bar\mu_0(\d x):= \ff{\e^{-V(x)-2\vv^{-1}\tt\rr (x)^2}\d x}{\int_{\R^d}\e^{-V(x)-2\vv^{-1}\tt\rr(x)^2}\d x},$$ which  is a probability measure by \eqref{PQ2''}.
Let $$\bar Z_0(x)= Z_0(x)- 2\vv^{-1} a(x)\nn \tt\rr(x)^2,\ \ \bar Z(x)= Z(x)+ 2\vv^{-1}a(x)\nn \tt\rr(x)^2.$$  By \eqref{PQ2''} we have
$\bar\mu_0(|\si^*\nn\psi|^2)<\infty$, so that $f_n:= (n- \psi)^+\land 1 \to 1$ in $L^2(\bar\mu_0)$ and
$$\lim_{n\to\infty}  \bar\mu_0(|\si^*\nn f_n|^2)=\lim_{n\to\infty} \int_{1+n\ge \psi\ge n} |\si^*\nn\psi|^2 \d\bar\mu_0 =0.$$ Thus, $1\in H_\si^{2,1}(\bar\mu_0)$ and $\bar\EE_0(1,1) =0.$ Then by
Lemma \ref{P1} for $(\bar Z_0, \bar Z, \bar\mu_0)$ in place of $(Z_0, Z, \mu_0)$, and due to \eqref{PQ2''},  it remains  to prove $\bar\mu_0(\e^{\vv'|\si^{-1}\bar Z|^2})<\infty$ for some $\vv'>0.$ Since
$|\si^*\nn \tt\rr|^2\le 2$,  we have
$$|\si^{-1}\bar Z|^2(x)\le 2 |\si^{-1}Z|^2(x)+ 8\vv^{-2}|(\si^*\nn \tt\rr(x)^2|^2(x)\le 2 |\si^{-1}Z|^2(x)+64\vv^{-2}\tt\rr(x)^2.$$ By \eqref{PQ2},
for $\vv'\in (0, \ff{\vv}{64}]$ we have
\beg{align*} \bar\mu_0(\e^{\vv' |\si^{-1}\bar Z|^2})&\le  \ff 1 {\int_{\R^d}\e^{-V(x)-2\vv^{-1}\tt\rr(x)^2}\d x} \int_{\R^d} \e^{2\vv' |\si^{-1}Z|^2(x) +  64\vv' \vv^{-2}  \tt\rr(x)^2 -  V(x) -2\vv^{-1}\tt\rr(x)^2}\d x\\
& \le  \ff 1 {\int_{\R^d}\e^{-V(x)-2\vv^{-1}\tt\rr(x)^2}\d x}\int_{\R^d} \e^{\vv |\si^{-1}Z|^2(x) - V(x)-\vv^{-1}\tt\rr(x)^2}\d x<\infty.\end{align*} Therefore,  the proof is finished.
\end{proof}

\beg{proof}[Proof of Corollary $\ref{C2.2}$] Let $n\ge 1$ such that $B_n:= \{|\cdot|\le n\}\supset D$. It suffices to show that for any
$l\ge n+1$ and any $x\in S_l:=\{|\cdot|=l\}$, the solution $\bar X_t^x$ to \eqref{E1} is non-explosive.  Let
$$\zeta^x= \lim_{m\to\infty}\inf\{t\ge 0: |\bar X_t^x|\ge m\},\ \ \si_n^x=\inf\{t\ge 0: |\bar X_t|\le n\},\ \ m>l\ge n+1, x\in S_l.$$
Let $\tt X_t^x$ solve the SDE \eqref{E1} for $Z1_{B_n^c}$ in place of $Z$. Due to \eqref{PQ*}, Theorem \ref{T1.0} applies to $\tt X_t$. In particular,  $\tt X_t^x$ is non-explosive, i.e.
\beq\label{NE}\tt\zeta^x:= \lim_{m\to\infty}\inf\{t\ge 0: |\tt X_t^x|\ge m\} =\infty,\end{equation} where and in the following,  $\inf\emptyset:=\infty$. Moreover, since $|Z|\in L_{loc}^p(\d x)$ for some $p>d$, \cite[Theorem 1.1]{Zhang2}  implies the pathwise uniqueness of the SDE \eqref{E1}. So,
$$\tt X_t^x= \bar X_t^x,\ \ t\le \si_n^x.$$ Then
\beq\label{YZ0} \si_n^x=\tt \si_n^x:=\inf\{t\ge 0: |\tt X_t^x|\le n\} \end{equation}  and
\beq\label{YZ1} \zeta^x=\tt \zeta^x \ \ \text{if}\ \ \zeta^x\le \si_n^x.\end{equation}
Obviously, for
$$\theta_n^x:=\inf\{t\ge \si_n^x:\ |X_t^x|\ge l\}  $$ we have
\beq\label{YZ2} \{\si_n^x< \zeta^x\}=\{\theta_n^x<\zeta^x\}.\end{equation}
By \eqref{NE}, \eqref{YZ2} and    the strong Markov property ensured by the uniqueness (see \cite[Theorem 5.1]{IW}),   we have
\beg{align*}&  \P(\zeta^x\le T) =   \P(\zeta^x\le T, \si_n^x\ge \zeta^x)+\P(\zeta^x\le T, \si_n^x<\zeta^x) \\
&\le   \P(\tt\zeta^x\le T)+ \P(\theta_n^x<\zeta^x\le T) = \E\big[1_{\{\theta_n^x\le T\}}\P(\theta_n^x<\zeta^x\le T|\F_{\theta_n^x})\big]\\
&= \E\big[1_{\{\theta_n^x\le T\}}\{\P(\zeta^z\le T-s)|_{s= \theta_n^x, z= X_{\theta_n^x}^x}\}\big]\\
&\le \P(\theta_n^x\le T) \sup_{z\in S_l} \P(\zeta^z\le T)\le \P(\si_n^x\le T) \sup_{z\in S_l} \P(\zeta^z\le T),\ \ T>0, x\in S_l.\end{align*} Combining this with \eqref{YZ0} we obtain
 \beq\label{YZ3} \sup_{x\in S_l} \P(\zeta^x\le T) \le \Big\{\sup_{x\in S_l}\P(\tt\si_n^x\le T)\Big\}\sup_{x\in S_l} \P(\zeta^x\le T),\ \ T>0. \end{equation}
Let $\tt P_t^\O$ be the Dirichlet semigroup of $\tt X_t^\cdot$ for $\O= B_n^c$. By applying  Theorem \ref{T1.0} for $Z1_{B_n^c}$ in place of $Z$, we obtain
$$\P(\tt\si_n^x\le T)= 1- \P(\tt\si_n^x>T)= 1- \tt P_T^\O 1(x)<1$$ and that $\P(\tt\si_n^x\le T)$ is continuous in $x\in \O$. So,
$$\vv_T:= \sup_{x\in S_l} \P(\si_n^x\le T) <1.$$
This together with \eqref{YZ3} implies  $\P(\zeta^x\le T)=0$ for any $T>0$ and $x\in S_l.$ Since $l\ge n+1$ is arbitrary and the   solution is continuous, we have $\P(\zeta^x=\infty)=1$ for all $x\in \R^d$.
\end{proof}

\section{Proofs of Theorem \ref{T1.1} and Theorem \ref{T1.1'}}

Since the uniqueness of invariant probability measure is ensured by the irreducibility and the strong Feller property, we only prove the existence and regularity estimates on the density. The  new technique in the proof of  the existence is to reduce the usual tightness condition to the boundedness of a Feyman-Kac semigroup, which follows from the hypercontractivity of $P_t^0$ under the given integrability condition. Moreover,
to estimate the derivative of the density, the  formula \eqref{Q2}  below will play a crucial role.

\beg{lem}\label{L2.2} Let $V\in W_{loc}^{1,1}(\d x)$ and $\si \in W_{loc}^{1,1}(\R^d\to \R^d\otimes\R^d; \d x)$ such that    $\mu_0(\d x)=\e^{-V(x)}\d x$ is a probability measure satisfying \eqref{SOB} and   the Poincar\'e inequality
\beq\label{P} \mu_0(f^2)\le C\mu_0(|\si^*\nn f|^2)+\mu_0(f)^2,\ \ f\in C_0^\infty(\R^d)\end{equation}   for some constant $C>0$. Let $L_0$ be in $\eqref{L0}$  and let $L:=L_0+Z\cdot\nn$ for some measurable $Z:\R^d\to\R^d$.
  If $Z$ has compact support and   $|Z|+|\nn \si|\in L^p(\d x)$ for some $p\in [2,\infty)\cap (d,\infty)$, then  any invariant probability measure $\mu$ of $L$ is absolutely continuous with respect to $\mu_0$ with density
$\rr:=\ff{\d\mu}{\d\mu_0}\in H^{2,1}_\si(\mu_0)$ satisfying
\beq\label{Q2*} \mu_0(\rr^2+ |\si^*\nn\rr|^2)\le (C+1)\, \mu_0(\rr^2|\si^*Z|^2)<\infty.\end{equation} Moreover,
\beq\label{Q2} \int_{\R^d} \<\si^*\nn f, \si^*\nn   \rr\> \d\mu_0= \int_{\R^d} \<Z, \nn  f\>\,\d\mu,\ \  f\in H^{2,1}_\si(\mu_0).\end{equation}
\end{lem}

\beg{proof}   Let $\mu$ be an invariant  probability measure of $L$.    Since  $|Z|+|\nn \si|$ is in $ L_{loc}^p(\d x)$
for some $p\in [2,\infty)\cap (d,\infty)$, by the local boundedness of $Z_0$ so is $|Z+Z_0|$. Then according to \cite[Corollary 1.2.8]{BKR}, for any   invariant probability measure $\mu$ of $L$,   $\mu(\d x)=\hat \rr(x)\d x$ holds for some $\hat\rr\in W_{loc}^{p,1}(\d x)$. Since $\mu_0(\d x)= \e^{-V(x)}\d x$ and $V\in C^1(\R^d)$, this implies
$\mu=\rr\mu_0$ for some $\rr\in W_{loc}^{2,1}(\d x).$ In particular, we may take a continuous version $\rr$  which is thus    locally bounded. By the integration by parts formula,
 \beq\label{EP} \beg{split} &\int_{\R^d} \<\si^*\nn \rr,\si^*\nn f\>\d\mu_0= -\int_{\R^d} \rr L_0 f\d\mu_0\\
&= -\int_{\R^d} Lf\d\mu +\int_{\R^d}\<Z, \nn f\>\,\d\mu= \int_{\R^d}\<\si^{-1}Z, \si^*\nn f\>\rr \,\d\mu_0,\ \ f\in C_0^\infty(\R^d).\end{split}\end{equation}
Since $Z$   has compact support with $|Z|\in L^2(\d x)$, and $\rr+\|\si^{-1}\|$ is locally bounded,  \eqref{EP}     implies
$$\bigg|\int_{\R^d} \<\si^*\nn \rr,\si^*\nn f\>\d\mu_0\bigg|\le \mu_0( \rr^2|\si^{-1}Z|^2)^{\ff 1 2} \mu_0(|\si^*\nn f|^2)^{\ff 1 2}<\infty,\ \ f\in C_0^\infty(\R^d).$$ Hence, $\mu_0(|\si^*\nn \rr|^2)\le \mu_0(\rr^2|\si^{-1}Z|^2)<\infty$. This and  \eqref{SOB} imply  $\rr\land N\in H^{2,1}_\si(\mu_0)$ for any $N\in (0,\infty)$. By the Poincar\'e inequality \eqref{P} we obtain
$$\mu_0((\rr\land N)^2) \le C\mu_0(|\si^*\nn (\rr\land N)|^2)+\mu_0(\rr)^2\le  C\mu_0(|\si^*\nn \rr|^2)+1<\infty,\ \ N\in (0,\infty).$$ By letting $N\to \infty$ we prove $\rr\in   H^{2,1}_\si(\mu_0)$ and \eqref{Q2*}, so that \eqref{Q2} follows from   \eqref{EP}.
\end{proof}

Below we will often  use the following version of Young's inequality  on a probability space
$(E,\B, \nu)$ (see \cite[Lemma 2.4]{ATW09}):
\beq\label{Yang} \nu(fg)\le \log \nu(\e^f) + \nu(g\log g),\ \ f,g\ge 0, \nu(g)=1.\end{equation}
The next lemma ensures the existence of invariant probability measure of $P_t$ for bounded $\si^{-1}Z$.

\beg{lem}\label{LNN} Assume {\bf (H1)}. If $\si^{-1}Z$ is bounded then the Markov semigroup $P_t$ associated to the SDE $\eqref{E1}$ has a unique invariant probability measure.
\end{lem}

\beg{proof} According to (b) in the proof of Lemma \ref{P1}, $P_t$ has at most one invariant probability measure. So, it suffices to prove the existence. Letting $\mu_0P_t$ be the distribution at time $t$ of the solution to \eqref{E1} with initial distribution $\mu_0$, we intend to show that the sequence $\{\ff 1 n\int_0^n \mu_0P_t\d t\}_{n\ge 1}$ is tight, so that the weak limit of a weakly convergent subsequence provides an invariant probability measure of $P_t$.
To this end, it suffices to find out a positive compact function $F$ on $\R^d$ such that
\beq\label{VV}\ff 1 n \int_0^n \mu_0(P_t F)\,\d t\le C,\ \ n\ge 1\end{equation} holds for some constant $C>0$.

According to Gross \cite{Gross}, {\bf (H1)} implies the hyperboundedness of $P_t^0$. Precisely, by \cite[Theorem 1]{Gross}   (see for instance also \cite[Theorem 5.1.4]{W05}),
we have
\beq\label{HPC} \|P_{t}^0\|_{L^q(\mu_0)\to L^{q(t)}(\mu_0)}\le \exp\Big[\bb \Big(\ff 1 q -\ff 1 {q(t)}\Big)\Big],\ \ t>0, q>1, q(t):=1+(q-1)\e^{\ff{4t}\kk}.\end{equation}   Since $\mu_0$ is a probability measure, there exists a  compact function $W\ge 1$ such that $\mu_0(W)<\infty$. Letting $F= \ss{\log W}$ which is again a compact function, we have $\mu_0(\e^{V^2})<\infty$.  We now prove \eqref{VV} for this function $F$. To this end, we consider the Feyman-Kac semigroup
$$P_t^Ff(x):= \E \Big[f(X_t^x) \e^{\int_0^tF(X_s^x)\d s}\Big],\ \ t\ge 0, x\in \R^d.$$ Since $\mu_0(\e^{F^2})<\infty$, $P_t^F$ is  a bounded linear operator from $L^p(\mu_0)$ to $L^1(\mu_0)$ for every $t\ge 0$ and $p>1$.    We first observe that $P_t^F$ is bounded on $L^p(\mu_0)$ for any $t\ge 0$ and $p>1$. Let $q=\ss p$. For any non-negative $f\in L^p(\mu_0)$, by Schwarz's and Jensen's inequalities, and     that $\mu_0$ is $P_t^0$-invariant,   we have

\beg{align*}  &\mu_0\big(|P_t^Ff|^p\big)  = \int_{\R^d} \Big(\E\Big[f(X_t^x)\e^{\int_0^t F(X_s^x)\d s}\Big]\Big)^p\mu_0(\d x)\\
&\le  \int_{\R^d} \bigg(\big\{\E f^{q}(X_t^x)\big\} \Big\{\E \e^{\ff{q}{q-1}\int_0^t F(X_s^x)\d s}\Big\}^{q-1} \bigg)^{q}\mu_0(\d x)\\
  &= \int_{\R^d}    (P_t^0 f^q)^{q} \bigg\{\ff 1 t \int_0^t P_s^0 \e^{\ff{qt}{q-1}F}\d s \bigg\}^{q(q-1)} \d\mu_0\\
 & \le  \big\{\mu_0((P_t^0f^q)^{q(t)})\big\}^{\ff q{q(t)}} \bigg\{\int_{\R^d}  \bigg(\ff 1 t \int_0^t P_s^0\e^{\ff{qt}{q-1}F} \d s\bigg)^{\ff {q(t)q(q-1)}{q(t)-q}}\d\mu_0\bigg\}^{\ff{q(t)-q}{q(t)}}\\
   &\le   \|P_{t}^0\|_{L^q(\mu_0)\to L^{q(t)}(\mu_0)}^q \mu_0\big(\{f^q\}^q\big) \max\bigg\{ \Big(\mu_0\Big(\e^{\ff{q^2q(t)t}{q(t)-q}F}\Big)\Big)^{\ff{q(t)-q}{q(t)}}, \ \Big(\mu_0\Big(\e^{\ff{qt}{q-1}F}\Big)\Big)^{q(q-1)}\bigg\}\\
   &= \|P_{t}^0\|_{L^q(\mu_0)\to L^{q(t)}(\mu_0)}^q \max\bigg\{ \Big(\mu_0\Big(\e^{\ff{q^2q(t)t}{q(t)-q}F}\Big)\Big)^{\ff{q(t)-q}{q(t)}}, \ \Big(\mu_0\Big(\e^{\ff{qt}{q-1}F}\Big)\Big)^{q(q-1)}\bigg\}\mu_0(f^p),\ \ t\ge 0.\end{align*}
By  $\mu_0(\e^{F^2})<\infty$ and \eqref{HPC}, this implies
$\|P_t^F\|_{L^p(\mu_0)}  <\infty $ for any $t>0,$ and moreover, $\limsup_{t\downarrow 0}\|P_t^F\|_{L^p(\mu_0)}\le 1.$ Since $F\ge 0$ implies $P_t^F1\ge 1$, we have $\lim_{t\downarrow 0}\|P_t^F\|_{L^p(\mu_0)}= 1.$ In particular, by taking $p=2$ and using the semigroup property, we obtain
\beq\label{BDD}\E \e^{\int_0^n F(X_t^{\mu_0})}= \mu_0(P_n^F 1)\le \|P_{n}^F\|_{L^2(\mu_0)}\le \|P_1^F \|_{L^2(\mu_0)}^n=:c_0^n<\infty,\ \ n\ge 1,\end{equation}
where $X_t^{\mu_0}$ is the solution to \eqref{E0} with initial distribution $\mu_0$. Now, define
$$R_n =\exp\bigg[\ff 1 {\ss 2}\int_0^n \<(\si^{-1}Z)(X_s^{\mu_0}),\d W_s\> -\ff 1 4\int_0^n |(\si^{-1}Z)(X_s^{\mu_0})|^2\d s\bigg],\ \ n\ge 0.$$ Since $\si^{-1}Z$ is bounded, by Girsanov's theorem we have
$$\mu_0(P_t F)= \E\big\{F(X_t^{\mu_0})R_n\big\},\ \ t\in [0,n].$$ Then   \eqref{Yang} and \eqref{BDD} imply
\beq\label{TH}  \beg{split}  &\ff 1 n \int_0^n \mu_0(P_t F)\d t=  \ff 1 n \int_0^n \E\big\{F(X_t^{\mu_0})R_n\big\} \d t\\
&\le \ff 1 n \log \E \e^{\int_0^n F(X_t^{\mu_0})\d t} + \ff 1 n \E\big\{R_n\log R_n\} \le c_0+ \ff 1 n \E\big\{R_n\log R_n\}.\end{split}\end{equation}
 Since by Girsanov's theorem
 $$\tt W_t:= W_t -\ff 1 {\ss 2}\int_0^t (\si^{-1}Z)(X_s^{\mu_0})\d s,\ \ t\in [0,n]$$ is a $d$-dimensional Brownian motion under the probability $\Q_n:= R_n\P$, we have
\beg{align*}&\E\big\{R_n\log R_n\}= \E_{\Q_n}\log R_n\\
&= \E_{\Q_n}\bigg(\ff 1 {\ss 2}\int_0^n \<(\si^{-1}Z)(X_s^{\mu_0}), \d \tt W_s\> +\ff 1 4 \int_0^n |(\si^{-1}Z)(X_s^{\mu_0})|^2\d s\bigg)\le \ff {n\|\si^{-1}Z\|_\infty^2} 4.\end{align*}
 Combining this with \eqref{TH}, we prove \eqref{VV}, and hence finish the proof. \end{proof}

\

\beg{proof}[Proof of Theorem $\ref{T1.1}(1)$] By Lemma \ref{P1}, {\bf (H1)} implies that \eqref{E1} has a unique  non-explosive solution and the associated Markov semigroup $P_t$ is strong Feller with at most one invariant probability measure. To apply Lemma $\ref{L2.2}$, we first consider bounded $Z$ with compact support, then pass to the general situation  by using  an approximation argument.

(a) Let $Z$ be bounded  with compact support.  By Lemma \ref{LNN},  $P_t$ has a unique  invariant probability measure $\mu$. In particular, $L^*\mu=0$, so that by  Lemma \ref{L2.2}(1) we have  $\mu=\rr\mu_0$ for some $\rr\in H^{2,1}_\si(\mu_0)$ such that \eqref{Q2} holds.

Since $\rr\in H^{2,1}_\si(\mu_0)$,  $f:=   \log (\rr +\dd)  \in H^{2,1}_\si(\mu_0)$ for all $\dd>0$. Taking this $f$ in \eqref{Q2} we obtain
\beg{equation*}\beg{split} &  \int_{\R^d} \ff{|\si^*\nn \rr|^2}{ \rr +\dd} \d\mu_0\le   \int_{\R^d} \{|\si^{-1}Z| \cdot  |\si^*\nn \log(\rr+\dd)|\} \,\d\mu\\
&=   \int_{\R^d} \{|\si^{-1}Z| \cdot |\si^*\nn \log(\rr+\dd)|\} \rr\,\d\mu_0  \le    \bigg(\int_{\R^d} \rr |\si^{-1}Z|^2 \d\mu_0\bigg)^{\ff 1 2}
 \bigg(\int_{\R^d}  \ff{\rr|\si^*\nn   \rr|^2}{( \rr +\dd)^2}\,\d\mu_0  \bigg)^{\ff 1 2}\\
  &\le  \bigg(\int_{\R^d} \rr |\si^{-1}Z|^2 \d\mu_0\bigg)^{\ff 1 2}
 \bigg(\int_{\R^d}  \ff{ |\si^*\nn   \rr|^2}{ \rr +\dd}\, \d\mu_0  \bigg)^{\ff 1 2},\ \ \dd>0.\end{split}\end{equation*}
 Since $\mu_0(\ff{|\si^*\nn   \rr|^2}{ \rr +\dd})<\infty$ due to $\rr\in H_\si^{2,1}(\mu_0)$, this implies
$$\int_{\R^d} \ff{|\si^*\nn   \rr|^2}{ \rr +\dd} \d\mu_0\le  \int_{\R^d} \rr |\si^{-1}Z|^2 \d\mu_0,\ \   \dd>0.$$ By letting $\dd\to 0$ we obtain
\beq\label{OP}\int_{\R^d} \big|\si^*\nn\ss{\rr}\big|^2 \d\mu_0\le \ff {1} {4}\int_{\R^d} \rr |\si^{-1} Z|^2\d\mu_0<\infty\end{equation} since $\si^{-1}Z$ is bounded and $\mu_0(\rr)=1.$ So,  $\ss\rr\in H^{2,1}_\si(\mu_0)$ by \eqref{SOB},  and
the log-Sobolev inequality \eqref{LS} implies
\beq\label{SLL} \mu(\rr\log\rr)\le \kk \int_{\R^d}  \big|\si^*\nn\ss\rr\big|^2\d\mu_0+\bb.\end{equation}
Combining this with \eqref{OP} and   the Young inequality \eqref{Yang},  we obtain
\beg{equation*}\beg{split}&\mu_0(|\si^*\nn\ss\rr|^2)\le \ff 1 {4  \ll} \log\mu_0(\e^{\ll |\si^{-1}Z|^2}) + \ff 1 {4  \ll} \mu_0(\rr\log\rr)\\
&\le \ff 1 {4 \ll} \log\mu_0(\e^{\ll |\si^{-1}Z|^2})+ \ff \kk {4 \ll} \mu_0(|\si^*\nn\ss\rr\big|^2)+ \ff{\bb}{4 \ll}.\end{split}\end{equation*} This and \eqref{OP} imply \eqref{PD}.

Similarly, $\rr\in H^{2,1}_\si(\mu_0)$ implies $f=    (\rr +\dd)^{-1}\in H^{2,1}_\si(\mu_0)$ for $\dd>0$, so that by \eqref{Q2} we have
\beg{equation*}\beg{split}   \int_{\R^d} \ff{|\si^*\nn  \rr|^2}{( \rr +\dd)^2} \d\mu_0
& \le    \int_{\R^d} \big\{\rr |\si^{-1}Z| \cdot |\si^*\nn (\rr+\dd)^{-1}|\big\}\d\mu_0\\
&\le   \bigg(\int_{\R^d} \ff{ (\rr |\si^{-1}Z|)^2}{(\rr+\dd)^2}\d\mu_0\bigg)^{\ff 1 2}
 \bigg(\int_{\R^d} \ff{|\si^*\nn  \rr|^2}{( \rr +\dd)^2} \d\mu_0  \bigg)^{\ff 1 2}\\
& \le \ss{\mu_0(|\si^{-1}Z|^2)}  \bigg(\int_{\R^d} \ff{|\si^*\nn  \rr|^2}{( \rr +\dd)^2} \d\mu_0  \bigg)^{\ff 1 2},\ \ \dd>0.\end{split}\end{equation*}
Therefore, \eqref{2.4} holds.

Finally, by \cite{[5]} the density function $\rr$ is strictly positive, so that by \eqref{2.4} and $H_\si^{2,1}(\mu_0)= W_\si^{2,1}(\mu_0)$ we have $\log \rr\in H_\si^{2,1}(\mu_0)$ if $\log\rr\in L^2(\mu_0)$. To prove $\mu_0(|\log \rr|^2)<\infty$, we use the Poincar\'e inequality. As explained above that the defective log-Sobolev inequality implies that the spectrum of $L_0$ is discrete, by the irreducibility of  the Dirichlet form we see that $L_0$ has a spectral gap, equivalently, the Poincar\'e inequality
$$\mu_0(f^2)\le C \mu_0(|\si^* \nn f|^2) +\mu(f)^2,\ \ f\in H_\si^{2,1}(\mu_0)$$ holds for some constant $C>0.$ Since
$\rr$ is strictly positive, we take $\vv\in (0,1)$ such that $\mu_0(\rr\le\vv)\le \ff 1 4.$ By \eqref{2.4} and $\mu_0(\rr)=1$, for any $\dd>0$  we have $\log(\rr+\dd)\in H_\si^{2,1}(\mu_0)$. Moreover, by     the Poincar\'e inequality, \eqref{2.4} and  $|\log (\rr+\dd)|\le \rr +\dd+\log\vv^{-1}$ for $\rr\ge\vv$, there exist constants $C_1,C_2>0$ such that
\beg{align*} &\mu_0(|\log (\rr+\dd)|^2)  \le C \mu_0(|\si^*\nn\log (\rr+\dd)|^2) +\mu_0(\log(\rr+\dd))^2\\
&\le C_1 + 2 \mu_0(\log(\rr+\dd)1_{\{\rr\le\vv\}})^2 + 2 \mu_0(\log(\rr+\dd)1_{\{\rr>\vv\}})^2\\
&\le C_1 + 2\mu_0(|\log(\rr+\dd)|^2) \mu_0(\rr\le\vv) + 2 \mu_0( \rr+\dd+\log \vv^{-1} )^2\\
&\le  \ff 1 2 \mu_0(|\log(\rr+\dd)|^2) +C_2,\ \ \dd\in (0,1).\end{align*}
Since $\mu(|\log(\rr+\dd)|^2)<\infty$, this implies
$$\mu(|\log\rr|^2)=\lim_{\dd\downarrow 0} \mu(|\log(\rr+\dd)|^2)\le 2 C_2 <\infty.$$

 (b) In general, for any $n\ge 1$ let
$$Z_n(x)=  1_{\{|x|+|Z(x)|\le n\}}Z(x),\ \ L_n= L_0+Z_n\cdot \nn.$$
By (a) and $|\si^{-1}Z_n|\le |\si^{-1}Z|$, $L_n$ has an invariant probability measure  $\d\mu_n =\rr_n\d\mu_0$ such that
\beg{align*} &\mu_0\big(|\si^*\nn\ss{\rr_n}|^2\big) \le \ff{1}{4  \ll- \kk} \big\{\log\mu_0(\e^{\ll|\si^{-1}Z|^2}) +\bb\big\}<\infty,\\
&\mu_0(|\si^*\nn \log \rr_n|^2) \le  \mu_0(|\si^{-1}Z|^2)<\infty.\end{align*}
Then the family $\{\ss{\rr_n}\}_{n\ge 1}$ is bounded in $H^{2,1}_\si (\mu_0)$.  Moreover, the defective log-Sobolev inequality \eqref{LS}
  implies the existence of a super Poincar\'e inequality, and hence   the essential spectrum of $L_0$ is empty, see
  \cite[Theorem  2.1 and Corollary 3.3]{W00a}. So,  $H^{2,1}_\si(\mu_0)$ is compactly embedded into $L^2(\mu_0)$, i.e.
   a bounded set in $H^{2,1}_\si(\mu_0)$ is relatively compact in $L^2(\mu_0)$. Therefore, for some subsequence $n_k\to\infty$ we have
   $\ss{\rr_{n_k}}\to \ss{\rr}$ in $L^2(\mu_0)$ for some nonnegative $\rr$ which satisfies \eqref{PD} and \eqref{2.4}. In particular, $\rr_{n_k}\to\rr$ in $L^1(\mu_0)$ so that $\mu:=\rr\mu_0$ is a probability measure. Moreover,   by using the Poincar\'e inequality as in (a), we prove $\log\rr\in L^2(\mu_0)$ so that $\log\rr \in H_\si^{2,1}(\mu_0)$. It remains to show that $L^*\mu=0.$

Since $(L_{n_k})^* \mu_{n_k}=0$, for any $f\in C_0^\infty(\R^d)$, there exists a constant $C>0$ and a compact set $D$ such that
\beq\label{WF} \beg{split}&\bigg|\int_{\R^d} Lf\d\mu\bigg| =   \bigg|\int_{\R^d} \big(\rr Lf-   \rr_{n_k}  L_{n_k}f \big)\d\mu_0\bigg|\\
&\le C    \int_{D} \Big\{|Z-Z_{n_k}|  \rr   +   (1+  |Z|) | \rr_{n_k}-\rr| \Big\}\d\mu_0. \end{split}\end{equation} Since $\mu_0(\e^{\ll |\si^{-1}Z|^2})<\infty$,
we have  $|Z_n|\le |Z|\in L_{loc}^q(\d x)$ for any $q>1$. Then
      $\mu_0(1_D|Z-Z_n|^q)\to 0$ as $n\to\infty$ holds for any $q>1$. Moreover,     the local Harnack inequality (see \cite[Corollary 1.2.11]{BKR}) implies that
  $\{\rr_{n_k},\rr\}_{k\ge 0}$  is uniformly bounded on the compact set $D$.
Combining these with  $ \mu_0(|\rr_{n_k}-\rr|)\to 0$, we may use   the dominated convergence theorem to prove $\mu(Lf)=0$
by taking $k\to \infty$ in \eqref{WF}. Therefore,    $L^*\mu=0$.
 Then the proof is complete.
\end{proof}

\beg{proof}[Proof of Theorem $\ref{T1.1}(2)$] By Theorem \ref{T1.0},     the SDE \eqref{E1} has a unique solution and the associated semigroup $P_t$ is strong Feller  having at most one invariant probability measure. So, it suffices to prove that the above constructed   probability measure $\mu$    is the unique  invariant probability measure of $L$ and  $P_t$. This can be done according to \cite{ST} and \cite{BKR} as follows.

Let $b_0= Z_0+a\nn\log\rr$ and $b= Z+Z_0$. Then $L= {\rm tr}(a\nn^2)+ b\cdot\nn$, and $\hat L_0:=  {\rm tr}(a\nn^2)+ b_0\cdot\nn$ is symmetric in $L^2(\mu)$. Obviously, {\bf (H1)} and \eqref{PP2} imply that conditions (1.1$'$)-(1.3$'$) and (1.4) in \cite{ST} hold for $U=\R^d$; that is, $a_{ij}\in W_{loc}^{2,1}(\d x)$, $a$ is   locally uniformly positive definite, and $b\in L_{loc}^2(\d x)$.   Moreover, by the Young inequality \eqref{Yang},  \eqref{PD},  \eqref{PP2},  \eqref{MU0} and \eqref{SLL},  for small enough $r>0$ we have \beg{align*} &\mu(\|a\|+|b-b_0|) \le \mu_0(\rr |Z|+ \|\si\|\cdot |\si^*\nn \rr|+\rr\|\si\|^2)\\
 &\le \ff 1 2 \mu_0(\rr(|\si^{-1}Z|^2+ 3\|\si\|^2)) + \mu_0(\|\si\|\cdot |\si^*\nn \rr|)\\
&\le \ff 1 {2r} \mu_0(\rr\log\rr)+\ff 1 {2r} \log \mu_0(\e^{r(|\si^{-1} Z|^2 +3\|\si\|^2)})  +2\ss{\mu_0(\rr\|\si\|^2) \mu_0(|\si^*\nn\ss{\rr}|^2)}\\
&\le \ff 1 {2r} \mu_0(\rr\log\rr)+\ff 1 {2r} \log \mu_0(\e^{r(|\si^{-1}Z|^2+ 3\|\si\|^2)}) \\
&\quad  +2\ss{\{\vv^{-1}\mu_0(\rr\log\rr)+\vv^{-1}\log\mu_0(\e^{\vv \|\si\|^2}) \}\mu_0(|\si^*\nn\ss{\rr}|^2)}<\infty.\end{align*}
Therefore, by \cite[Theorem 1.5, Proposition 1.9 and Proposition 1.10(a)]{ST}, $(L, C_0^\infty(\R^d))$ has a unique closed extension in $L^1(\mu)$ which  generates a Markov $C_0$-semigroup
$T_t^\mu$ in $L^1(\mu)$ such that $\mu$ is an invariant probability measure. Then, according to \cite[Corollary 1.7.3]{BKR}, $\mu$ is the unique invariant probability measure of $L$.

 On the other hand, according to \cite[Theorem 3.5]{ST}, there is a standard Markov process $\{\bar\P_x\}_{x\in\R^d\cup\{\pp\}}$ which is continuous and  non-explosive for $\mu$-a.e. $x$, such that the associated semigroup $\bar P_t$ satisfies
  $$\int_0^\infty \e^{-\ll t} \bar P_tf\d t= \int_0^t \e^{-\ll t} T_t^\mu f\d t,\ \ \mu\text{-a.e.}$$ holds for   any $f\in \B_b(\R^d)$ and $\ll>0$.
    So,  for any $f\in \B_b(\R^d), \bar P_t f= T_t^\mu f$ holds $\d t\times\mu$-a.e. By the continuity of the process and the strong continuity of $T_t^\mu$ in $L^1(\mu)$, $\bar P_t f= T_t^\mu f \ \mu$-a.e. for any $t\ge 0$ and $f\in  C_b(\R^d)$, and hence also for $f\in L^1(\mu)$ since $C_b(\R^d)$ is dense in $L^1(\mu)$. That is, $\bar P_t$ is a $\mu$-version of $T_t^\mu$. In particular, $\mu$ is $\bar P_t$-invariant and the probability measure
    $$\bar\P_\mu:= \int_{\R^d} \P_x\mu(\d x)\ \text{on}\ \ \bar\OO:=C([0,\infty)\to\R^d)$$ solves the martingale problem of $(L, C_0^\infty(\R^d))$, so that under this probability space the coordinate process $\bar X_t(\bar\oo):= \bar\oo_t$ for $t\ge 0$ and $\bar\oo\in \bar\OO$ is a weak solution to \eqref{E1}
    with initial distribution $\mu$ (c.f. \cite[Proposition 2.1]{IW} or \cite[\S 5.0]{SV}). By the uniqueness of solutions, this implies $\mu(P_tf)=\mu(\bar P_tf)$ for $t\ge 0$ and $f\in \B_b(\R^d)$. Therefore,     $\mu$ is an invariant probability measure of $P_t$.
  \end{proof}

\beg{proof}[Proof of Theorem \ref{T1.1'}] Obviously, the proof of Theorem \ref{T1.1}(2) also works if we replace {\bf (H1)} by {\bf (H$1'$)}. So, we only need to prove assertion (1).  Next,  by repeating  (b) in the proof of Theorem \ref{T1.1}(1), we may and do assume that $Z$ is bounded having compact support, and only prove that $L$ has an invariant probability  measure $\d\mu =\rr\d\mu_0$ with $\rr$ satisfying the required estimates
\eqref{*D1} and \eqref{*D2}. Here, the only thing we need to clarify is that in the right hand side of \eqref{WF} the term $(1+|Z|)$ should be replaced by $(1+|Z|+|\nn\si|)$ since $\nn \si$ is no longer locally bounded. This does not make any trouble since $|\nn\si|\in L_{loc}^2(\d x)$ by {\bf (H$1'$)}, and  $(\rr_{n_k}-\rr)1_D$ is uniformly bounded according to \cite[Corollarty 1.2.11]{BKR}.

Now, we assume that $Z$ is bounded with compact support. Let $\tt V\in C^\infty(\R^d)$ with $\|\tt V-V\|_\infty\le 1$, and let $\tt P_t$ be the Markov semigroup generated by $\DD-\nn \tt V$. Then $H^{2,1}(\mu_0)=H^{2,1}(\e^{-\tt V(x)}\d x),$ so that   {\bf (H$1'$)} together with the smoothness and positivity-preserving of $\tt P_t$ implies
\beq\label{**D} \beg{split} &\tt a_n:= \tt P_{\ff 1 n} a\in C^2(\R^d\to\R^d\otimes \R^d),\ \   \tt a_n\ge \aa I,\\
&\text{and}\   (\tt a_n)_{ij}\to a_{ij}\ \text{in}\ H^{2,1}(\mu_0)\cap L^{2p}(\mu_0), \ 1\le i,j\le d.\end{split}\end{equation}

Let $\tt L_n$ be defined as $L$ for $\tt a_n$ in place of $a$; that is,
$$\tt L_n= {\rm tr}(\tt a_n\nn^2) +\sum_{i,j=1}^d \big\{Z_i+ \pp_j (\tt a_n)_{ij} -(\tt a_n)_{ij} \pp_j V\big\} e_i.$$
  By Lemmas \ref{L2.2} and \ref{LNN}, $\tt L_n$ has an invariant probability measure $\tt \mu_n(\d x):= \tt\rr_n(\d x)\mu_0(\d x)$ with $\tt\rr_n\in H^{2,1}(\mu_0)$ such that
$$\mu_0(\tt\rr_n^2+ |\nn \tt\rr_n|^2) \le C\mu_0(\tt\rr_n^2 |Z|^2)<\infty.$$ According to \cite[Corollarty 1.2.11]{BKR}, $\{\tt\rr_n\}_{n\ge 1}$ is uniformly bounded on the compact set $D:=\text{supp}\,Z$, so   this implies that $\{\tt\rr_n\}_{n\ge 1}$ is bounded in $H^{2,1}(\mu_0)$, and hence $\tt\rr_{n_k}\to\rr$ in $L^2(\mu_0)$ for some subsequence $n_k\to\infty$ and some $\rr\in H^{2,1}(\mu_0)$. In particular, $\mu(\d x):=\rr(x)\d x$ is a probability measure. We intend to prove $L^*\mu=0$.

For any $f\in C_0^\infty(\R^d)$ there exists a constant $C(f)>0$ such that
\beg{align*}&|\tt L_n f-Lf|\le C(f) \big(\|\nn \tt a_n-\nn a\|+|\nn V|\cdot\|\tt a_n-a\|\big),\\
 &|\tt L_n f|\le C(f) (\|\nn \tt a_n\| + \|a\|\cdot |\nn V|).\end{align*}
By \eqref{**D}, $|\nn V|\in L^{\ff{2p}{p-1}}(\mu_0)$ included in {\bf (H$1'$)},  $\tt\rr_{n_k}\to \rr$ in $L^2(\mu_0)$, and $\tt L_n^*\tt\mu_n=0$, we are able to use the dominated convergence theorem to derive
$$|\mu(Lf)|= \lim_{k\to\infty}|\mu(Lf)-\tt\mu_{n_k}(\tt L_{n_k}f)|\le \limsup_{k\to\infty} \mu_0(|Lf-\tt L_{n_k}f|\rr+ |\tt L_{n_k} f|\cdot|\tt\rr_{n_k}-\rr|)=0.$$
So, $L^*\mu=0$.

Since \eqref{LS'} and $\tt a_n\ge \aa I$ imply \eqref{LS} for $\big(\ss{\tt a_n}, \ff{\kk'}\aa\big)$ in place of $(\si, \kk)$,  by Theorem \ref{T1.1} we have
\beg{align*} &\aa \mu_0\Big(|\nn\ss{\tt\rr_{n_k}}|^2\Big)\le \mu_0\Big(\Big|\ss{\tt a_{n_k}}\nn\ss{\tt\rr_{n_k}}\Big|^2\Big)\\
 &\qquad \le \ff{1}{4 \aa\ll- \ff{\kk'}\aa} \big\{\log\mu_0\big(\e^{\aa\ll|(\tt a_{n_k})^{-1/2}Z|^2}\big) +\bb\big\}\le \ff{\aa}{4\aa^2\ll-\kk'}\big\{\log\mu_0\big(\e^{\ll |Z|^2}\big) +\bb\big\},\\
& \aa \mu_0(|\nn\log \tt\rr_{n_k}|^2) \le \mu_0\big(\big|(\tt a_{n_k})^{-1/2}\nn\log\tt\rr_{n_k}\big|^2\big) \le  \ff 1 {\aa}\mu_0(|Z|^2).  \end{align*}
By   using $\rr_{n_k}+\dd$ to replace $\rr_{n_k}$, and  letting first $k\to\infty$ then $\dd\downarrow 0$,   we prove \eqref{*D1} and \eqref{*D2} from these two inequalities respectively.
\end{proof}

\section{Proofs of Theorem \ref{T4.1} and Theorem \ref{T4.1'}}

The following Sobolev embedding theorem is crucial in the proof. This result can be deduced from existing ones, for instance, \cite[Corollary 1.4]{RW} in the framework of generalized Mehler semigroup. We include below a brief proof by using the dimension-free Harnack inequality for   the O-U semigroup.

\beg{lem}\label{L4.1} Let $\eqref{EG}$ hold. Then $H^{2,1}(\mu_0)$ is compactly embedded into $L^2(\mu_0)$; i.e. bounded sets in $H^{2,1}(\mu_0)$ are relatively compact in $L^2(\mu_0).$ \end{lem}

\beg{proof}   Consider the linear  SPDE
 \beq\label{EE} \d X_t = -AX_t\d t +\ss 2\,\d W_t, \end{equation}
By $\eqref{EG}$, for any initial point $x$  this equation
has a unique mild solution
$$X_t^x= \e^{-At} x+ \ss 2 \int_0^t\e^{-A(t-s)}\d W_s,\ \ t\ge 0, $$ and the associated Markov semigroup
$$P_t^0 f(x):= \E f(X_t^x),\ \ t\ge 0, f\in \B_b(\H), x\in\H$$ is symmetric in $L^2(\mu_0)$ with Dirichlet form
$$\EE_0(f,g):=\mu_0(\<\nn f,\nn g\>),\ \ \ f,g\in H^{2,1}(\mu_0),$$ see for instance \cite{DZ}.   So, by the spectral theory, $H^{2,1}(\mu_0)$ is compactly embedded into $L^2(\mu_0)$ if and only if $P_t^0$ is compact for some (equivalently, all) $t>0$, both are equivalent to the absence of the essential spectrum of the generator.   By   \cite[Theorem 3.2.1]{Wbook2} with $b=0$ and $\si=\ss 2$ so that $K=0$ and $\ll=\ff 1 2$, $P_t^0$ satisfies the Harnack inequality
\beq\label{HN}(P_t^0f(x))^2\le (P_t^0f(y))^2 \e^{2|x-y|^2/t},\ \ t>0, x,y\in \H, f\in \B_b(\H),\end{equation}  which implies that $P_t^0$ has a density with respect to the invariant probability measure $\mu_0$ (see \cite[Theorem 1.4.1]{Wbook2}). Next,  it is well known that the Gaussian measure $\mu_0$ satisfies the log-Sobolev inequality (see for instance \cite{Gross})
\beq\label{LS'} \mu_0(f^2\log f^2)\le \ff 2 {\ll_1} \mu_0(|\nn f|^2),\ \ f\in H^{2,1}(\mu_0), \mu_0(f^2)=1.\end{equation}
This,  together with the existence of density of $P_t^0$  with respect to $\mu_0$ for any $t>0$, implies that  $P_t^0$ is compact in $L^2(\mu_0)$ for all $t>0$, see \cite[Theorem 1.2]{GW},   \cite[Theorems 1.1 and 3.1]{W00b} or \cite[Theorem 1.6.1]{Wbook}.\end{proof}

\beg{proof}[Proof of Theorem $\ref{T4.1}(1)$]  For any $n\ge 1$,   let $\H_{\<n\>}=\{x\in\H: \<x,e_i\>=0, 1\le i\le n\}$ be the  orthogonal complement of $\H_n:= {\rm span} \{e_1,\cdots, e_n\}$.  Let  $\pi_n: \H\to\H_n$ and $\pi_{\<n\>}: \H\to \H_{\<n\>}$ be orthogonal projections.     For convenience, besides the orthogonal decomposition  $\H= \H_n \oplus \H_{\<n\>}$ we may  regard $\H$ as the product space  $\H=\H_n\times \H_{\<n\>}$, so that $\mu_0= \mu_0^{(n)}\times \mu_0^{\<n\>}$ for  $\mu^{(n)}_0 =\mu_0\circ \pi_n^{-1}$ and $\mu^{\<n\>}_0
=\mu_0\circ \pi_{\<n\>}^{-1}$ being the marginal distributions of $\mu_0$ on $\H_n$ and $\H_{\<n\>}$ respectively.  Let
\beq\label{ZN}a_n(x)=\pi_n a(x),  \ \  Z_n(x)=\pi_n \int_{\H_{\<n\>}}Z(x,y)\mu_0^{\<n\>}(\d y),\ \ x\in\H_n. \end{equation} By
 {\bf (H2)} we have
 \beq\label{SI'} \<a_n v,v\>\ge \aa |v|^2,\ \ v\in \H_n,\end{equation} and due to Jensen's inequality,
 \beq\label{LB} \mu_0^{(n)}(\e^{\ll|Z_n|^2})\le\int_{\H_n} \e^{\ll \int_{\H_{\<n\>}} |Z(x,y)|^2\mu_0^{\<n\>}(\d y)}\mu_0^{(n)}(\d x) \le \int_\H \e^{\ll |Z|^2} \d \mu_0 <\infty,\ \ n\ge 1.\end{equation}
 Let $V_n(x)= \ff 1 2\sum_{i=1}^n \ll_i x_i^2$ and $L^{(n)}= L_0^{(n)}+Z_n\cdot\nn$ on $\H_n$, where
 $$L_0^{(n)} =   \sum_{i,j=1}^n \Big(a_{ij} \pp_i\pp_j +\big\{\pp_j a_{ij} -a_{ij}\pp_j V_n \big\}\pp_i\Big).$$   Noting that    \eqref{LS'} and {\bf (H2)} imply
\beq\label{LS"} \mu_0(f^2\log f^2)\le \ff 2 {\ll_1\aa} \mu_0(\<a\nn f,\nn f\>),\ \ f\in H^{2,1}(\mu_0), \mu_0(f^2)=1,\end{equation}
 and  that \eqref{SI'} implies $\aa|a_n^{- 1/2}Z_n|^2\le  |Z_n|^2,$
 we may apply Theorem \ref{T1.1}(1)   to $L_n$ on $\R^n\equiv \H_n$ for  $\kk= \ff 2 {\ll_1\aa},\bb=0$ and $\ll\aa$ in place of $\ll$, to conclude that  $L^{(n)}$ has an invariant  probability measure $\mu_n$ with density function $\rr_n:=\ff{\d\mu^{(n)}}{\d\mu_0^{(n)}}$ satisfying
 $\ss\rr_n\in H^{2,1}(\mu_0^{(n)})$ and
\beq \label{LB2}\beg{split}&\mu_0^{(n)}\big(\big|\nn \ss\rr_n\big|^2\big)\le \ff 1 {\aa} \mu_0^{(n)}(|\ss{a_n} \nn \ss{\rr_n}|^2)
  \le \ff{\ll_1}{4\aa^2\ll_1\ll- 2} \log\mu_0^{(n)}(\e^{\ll\aa |a_n^{-1/2}Z_n|^2})\\
  &\le\ff{\ll_1}{4\aa^2\ll_1\ll- 2} \log\mu_0^{(n)}(\e^{\ll|Z_n|^2})\le \ff{\ll_1}{4\aa^2\ll_1\ll- 2} \log\mu_0 (\e^{\ll|Z|^2}) <\infty,\ n\ge 1,\end{split}\end{equation}where the last step is due to Jensen's inequality and the definitions of $Z_n$ and $\mu_0^{(n)}$. Moreover,
\beq\label{LBB}\beg{split} &\mu_0^{(n)}(|\nn \log \rr_n|^2)\le \ff 1 \aa \mu_0^{(n)}(|\ss{a_n}\nn \log \rr_n|^2) \\
&\qquad \le  \ff{\mu_0^{(n)}(|a_n^{-  1/2}Z_n|^2)}{\aa}\le \ff{\mu_0^{(n)}(|Z_n|^2)}{\aa^2}\le \ff{\mu_0(|Z|^2)}{\aa^2}<\infty,\ \ n\ge 1. \end{split}\end{equation}
Letting $\bar\rr_n=\rr_n\circ\pi_n$, \eqref{LB2} implies that $\{\ss{\bar\rr_n}\}_{n\ge 1}$ is bounded in $H^{2,1}(\mu_0)$. By   Lemma \ref{L4.1},  there exists a subsequence $n_k\to\infty$ and some positive $\rr\in L^1(\mu_0)$ with $\ss\rr\in H^{2,1}(\mu_0)$ such that  $\ss{\bar\rr_{n_k}}\to \ss{\rr}$ in $L^2(\mu_0)$,    \eqref{PD'} and \eqref{2.4'} hold. Then  $\log\rr \in H^{2,1}(\mu_0)$ as shown in the end of the proof of Theorem \ref{T1.1}(1) using the Poincar\'e inequality.
In particular, $\mu:=\rr\mu_0$ is a probability measure on $\H$. It remains to show that $L^*\mu=0$.

By the definition of $Z_n$, we have   $\bar Z_n:= Z_n\circ\pi_n= \pi_n\mu_0(Z|\pi_n),$ where $\mu_0(\cdot|\pi_n)$ is the conditional expectation of $\mu_0$ given $\pi_n$. Since $\mu_0(|Z|^2)<\infty$, by the martingale converges theorem, $\mu_0(Z|\pi_n)\to Z$ in $L^2(\mu_0)$, and hence,   $\bar Z_n\to Z$ in $L^2(\mu_0)$ as well.    By the continuity of $a$,    $\bar a_n:= a_n \circ\pi_n \to a$ pointwise.   Noting that for any $f\in \F C_0^\infty$  there exist $l\in \mathbb N$ and a constant $C(f)>0$ such that
\beg{align*}& |\mu(Lf)|= |\mu(Lf)-\mu_{n_k}(L_{n_k}f)|\\
&\le C(f) \mu_0\Big(\rr\{|Z-\bar Z_{n_k} |+\sum_{i,j=1}^l |(a-\bar a_{n_k})_{ij}|\}\Big)+ C(f)\mu_0\Big(\Big\{|\bar Z_{n_k}|+\sum_{i,j=1}^l |(\bar a_{n_k})_{ij}|\Big\} |\rr-\bar\rr_{n_k}|\Big) \end{align*} holds for $n_k\ge l$, to prove $\mu(Lf)=0$ by using the dominated convergence theorem,    it suffices to
verify the uniform integrability of $\{\bar\rr_n( |\bar Z_n|+  |a_{ij}\circ\pi_n|)\}_{n\ge 1}$ in $L^1(\mu_0)$ for every $i,j\ge 1$. Obviously, for any $\vv\in (0,1)$ there exists a constant $C(\vv)>0$ such   that
$$(|\bar Z_n|+|a_{ij}\circ\pi_n|) \bar\rr_n
 \le \e^{\vv  |\bar Z_n|^{1+\vv}} +\e^{\vv|a_{ij}\circ\pi_n|^{1+\vv}}  + C \bar\rr_n\{\log (\e+ \bar\rr_n)\}^{\ff 1{1+\vv}},\ \ n\ge 1.$$ Since $\mu_0^{(n)}(f)= \mu_0(f\circ\pi_n)$ for $f\in L^1(\mu_0^{(n)})$, this implies the desired
the uniform integrability by \eqref{SII},  \eqref{LB},   \eqref{LB2}   and
$$\mu_0(\bar\rr_n\log \bar\rr_n)\le \ff 2 {\ll_1} \mu_0(|\nn \ss {\bar\rr_n}|^2)= \ff 2 {\ll_1} \mu_0(|\nn \ss {\rr_n}|^2)$$ due to the log-Sobolev inequality \eqref{LS'}. \end{proof}

\beg{proof}[Proof of Theorem $\ref{T4.1}(2)$] The desired assertion can be deduced  from \cite{ST}. Since $a$ is bounded and {\bf (H2)} holds, we have $H^{2,1}(\mu_0)=H^{2,1}_\si(\mu_0).$ Let $\mu$ be  a probability measure $\mu$ on $\H$ such that  the form
  $$\EE^\mu(f,g):= \mu(\<a\nn f,\nn g\>),\ \ \ f,g\in \F C_0^\infty $$ is closable in $L^2(\mu)$, and let   $(L^\mu,\D(L^\mu))$ be the generator of the closure $(\EE^\mu, H^{2,1}(\mu))$. Moreover,
  let $\bb\in L^2(\H\to\H;\mu)$ such that
  \beq\label{GSS20} \mu(\<\bb, \nn f\>)=0,\ \ f\in H^{2,1}(\mu).\end{equation}
  Then,  according to   Proposition 1.3, Theorem 1.9 and Proposition 1.10 in \cite[Part II]{ST}, we have the following assertions for $L:=L^\mu+\bb\cdot\nn$:
\beg{enumerate} \item[(i)] $(L, \F C_b^\infty)$ is dissipative and hence closable in $L^1(\mu)$, whose closure $(\bar L, \D(\bar L))$ generates a Markovian $C_0$-semigroup of contraction operators $(T_t)_{t\ge 0}$ on $L^1(\mu)$,   $\D(\bar L)\subset H^{2,1}(\mu)$, and
\beq\label{*Y20} \mu(\<\nn f, \bb-a\nn g\>) =\mu(g\bar L f),\ \ f\in\D(\bar L)\cap \B_b(\H), g\in H^{2,1}(\mu)\cap \B_b(\H).\end{equation}
\item[(ii)] There exists a standard continuous  Markov process $\{\bar\P_x \}_{x\in \H}$ whose semigroup $\bar P_t$ satisfies
\beq\label{II} \int_0^\infty \e^{-\ll t} \bar P_t f\d t= \int_0^\infty T_tf\d t,\ \ \mu\text{-a.e.},\ \ll>0, f\in\B_b(\H).\end{equation}   \end{enumerate} As shown in the proof of Theorem \ref{T1.1}(2),
\eqref{II}  implies  that  $\bar P_t$ is a $\mu$-version of $T_t$.

\

Now, let $L=L_0+Z\cdot\nn$ and $\mu=\rr\mu_0$ be in Theorem \ref{T4.1}. We intend to verify the above conditions such that assertions (i) and (ii) hold.

Firstly,  $\ss\rr\in H^{2,1}(\mu_0)$ implies  $\nn\log \rr\in L^2(\mu)$ and
  $$\mu_0(|\nn\rr|)\le 2\ss{\mu_0(|\nn \ss\rr|^2)\mu_0(\rr)} <\infty.$$
 Consider the operator
$$L^\mu:= L_0 +a\nn\log\rr,\ \ f\in \F C_0^\infty.$$ By the symmetry of $L_0$ in $L^2(\mu_0)$, the boundedness of $a$, $\nn\log\rr\in L^2(\mu_0), \nn\rr\in L^1(\mu_0)$ and noting that $H^{2,1}(\mu_0)$ is dense in $H^{1,1}(\mu_0)$,  we obtain
\beg{align*} &\mu(f L^\mu g)=\mu(f\<\nn\log\rr,a\nn g\>)+\mu_0(f\rr L_0g)\\
&=\mu(f\<\nn\log\rr,a\nn g\>) -\mu_0(\nn(f\rr),a\nn g\>)= -\mu(\<\nn f,a\nn g\>),\ \ f,g\in \F C_0^\infty. \end{align*}
 Thus, the form $(\EE^\mu, \F C_0^\infty)$   is closable in $L^2(\mu)$ with generator extending $(L^\mu,\F C_0^\infty)$.

 Next, let $\bb= Z-a\nn\log \rr.$ We have $L= L^\mu +\bb\cdot\nn$ on $\F C_0^\infty$.
Since $L^*\mu=0$ and $\mu_0(\<\nn\rr,\nn f\>)=-\mu_0(\rr L_0f)$ for $f\in \F C_0^\infty$, we have
\beq\label{GSS0}\beg{split}\mu(\<\bb, \nn f\>)&= \mu_0(\<\rr Z-a\nn\rr, \nn f\>) \\
&= \mu(\<Z,\nn f\>) + \mu_0(\rr L_0 f) =\mu(Lf)=0,\ \ f\in \F C_0^\infty.\end{split}\end{equation}
Noting that \eqref{PD'} and the boundedness of $a$ imply
 $$\mu(|a\nn\log\rr|^2)\le 4\|a\|^2 \mu_0(|\nn \ss \rr|^2) <\infty,$$ while by
 the Young inequality \eqref{Yang} and the log-Sobolev inequality \eqref{LS'}
\beg{align*} &\mu(|Z|^2)= \mu_0(\rr |Z|^2)\le \ff 1 \ll \log \mu_0(\e^{\ll|Z|^2}) + \ff 1 \ll \mu_0(\rr\log\rr)\\
&\le \ff 1 \ll \log \mu_0(\e^{\ll|Z|^2}) + \ff 2{\ll_1  \ll} \mu_0(|\nn\ss\rr|^2)+\ff 1 \ll<\infty,\end{align*} we have   $\mu(|\bb|^2)<\infty$
for $\bb:= Z-a\nn\log\rr$. So,
\eqref{GSS20}  follows from \eqref{GSS0}.

In conclusion, the above assertions (i) and (ii) hold for the present situation.
Combining    \eqref{GSS20} with \eqref{*Y20} for $g=1$ and $T_tf$ in place of $f$,   we obtain
$$ \ff{\d}{\d t}\mu( T_t f)= \mu(L  T_t f) =\mu(\<\nn T_t f,\bb\>)=0,\ \ f\in \F C_0^\infty, t\ge 0.$$ Therefore, $\mu$ is an invariant probability measure of $T_t,$  and the proof is finished since
 $\bar P_t$ is a $\mu$-version of $T_t$.
\end{proof}

\beg{proof}[Proof of Theorem $\ref{T4.1'}$] Since $V_n(x):= \ff 1 2\sum_{i=1}^n \ll_i x_i^2$ on $\H_n$ satisfies $|\nn V_n|\in L^1(\mu_0^{(n)})$ for all $q>1$,  {\bf (H$2'$)} and \eqref{LS"} imply that {\bf (H$1'$)} holds for $(a_n, V_n, \mu_0^{(n)})$ in place of $(a, V, \mu_0)$ with $\kk'= \ff 2 {\aa\ll_1}$ and $\bb=0.$ So, by  repeating the proof of Theorem \ref{T4.1} using Theorem \ref{T1.1'}  in place of Theorem \ref{T1.1}(1), we prove the desired assertions.
\end{proof}

\section{Proof of Theorem \ref{T4.2}}

We first prove the non-explosion of the weak solution constructed from the Girsnaov transform of the linear SPDE \eqref{EE}, then prove the strong Feller property of the associated Markov semigroup. The Feller property, together with   the pathwise uniqueness for $\mu_0$-a.e. starting points due to \cite{DFRV}, implies that the constructed Markov process is the unique Feller process   solving \eqref{E1} weakly.
Noting that in the present case we have $d=\infty$, the estimate \eqref{EST2} derived in the finite-dimensional case does not make sense. To construct the desired weak solution we need to establish a reasonable infinite-dimensional version of \eqref{EST2}. We will soon find out that this is non-trivial at all. If we start from the Harnack inequality $\eqref{HN}$,    it is  standard  that
$$(P_t f(x))^p \le \ff{\mu_0(f^p)}{\mu_0(\e^{-|x-\cdot|^p/t})}\approx \e^{c(x)/t}$$ for some constant $c(x)>0$ and small $t>0$. The hard point is that $\int_0^t\e^{c(x)/(ps)}\d s=\infty$  for any $t>0$ and $p>1$, so that the argument we used in the finite-dimensional case is invalid. To kill this high singularity  for small time $t$, we will use a refined version of the Harnack inequality and make a clever choice
of reference measure $\nu_t$ on $[0,t]$ to replace the Lebesgue measure.

 \subsection{Construction of the weak solution  }

 We first construct   weak solutions to \eqref{E1} using the Girsanov transform. For any $x\in\H$, let $X_t^x$ solve \eqref{EE} with $X_0=x$. Let
\beq\label{RX} R_{s,t}^x:=\exp\bigg[\ff 1 {\ss 2}\int_s^t\<Z(X_r^x), \d W_r\>-\ff 1 4 \int_s^t |Z(X_r^x)|^2\d r\bigg],\ \ t\ge s\ge 0.\end{equation} By Girsanov's theorem, if
 $(R_t^x)_{t\ge 0}:= (R_{0,t}^x)_{t\ge 0}$ is a martingale, then for any $T>0$ the process
 $$\widetilde W_t^x:= W_t -\ff 1 {\ss 2}\int_0^t Z(X_s^x)\d s,\ \ t\in [0,T]$$ is a cylindrical Brownian motion under the weighted probability $Q_T^x:= R_T^x\P$, so that   $(X_t^x, \widetilde W_t^x)_{t\in [0,T]}$ is a weak solution to \eqref{E2} starting at $x$. To prove that $(R_t^x)_{t\ge 0}$ is a martingale, it suffices to verify   the Novikov condition
\beq\label{AB0}\E\e^{\ff 1 4 \int_0^{t_0} |Z(X_s^x)|^2\d s} <\infty,\ \ x\in\H\end{equation}for some $t_0>0.$ Indeed, by the Markov property,   this condition   implies that $(R_{s,t}^x)_{t\in [s,s+t_0]}$ is a martingale for all $x\in\H$ and $s\ge 0$, and thus $(R_t^x)_{t\ge 0}$ is a martingale for all $x\in\H$ by induction:    if $(R_t^x)_{t\in [0, nt_0]}$ is a martingale for some $n\ge 1$, then for any $nt_0\le s<t\le (n+1)t_0$ we have
$$\E (R_t^x|\F_s)= R_s^x\E(R_{s,t}^x|\F_s)= R_s^x.$$
Therefore, the condition \eqref{AB0} implies that   $(X_t^x, \widetilde W_t^x)_{t\in [0,T]}$ is a weak solution to \eqref{E2} for any $T>0$ and $x\in\H$.  Let $ P_t(x,\d y)$ be the distribution of $X_t^x$ under $\Q_t^x$, and let
 \beq\label{GR} P_t f(x)= \E_{\Q_t^x} f(X_t^x)=\E \big\{f(X_t^x)R_t^x\big\},\ \ f\in \B_b(\H), t\ge 0, x\in\H.\end{equation} By the Markov property of $X_t$ under $\P$, it is easy to see that $P_t$ is a Markov semigroup on $\B_b(\H)$, i.e. $\{P_t(x,\d y): t\ge 0, x\in\H\}$ is a Markov transition kernel.

 To verify condition \eqref{AB0}, we introduce a refined version of the Harnack inequality \eqref{HN}. For each $i\ge 1$ let $P_t^{0,i}$ be the diffusion semigroup generated by $L_{0,i}f:= f''-\ll_i f'$ on $\R$. By \cite[Lemma 2.1]{W97} for $K=-\ll_i$ and $g(s)= \e^{-Ks}$, we have
 $$(P_t^{0,i}f(x))^p\le (P_t^{0,i}f^p(y))\exp\bigg[\ff{p\ll_i|x-y|^2}{2(p-1)(\e^{2\ll_it}-1)}\bigg],\ \ t>0, p>1, f\in \B^+(\R), x,y\in\R.$$
 By regarding $P_t^{0,i}$ as a linear operator on $\B_b(\H)$   acting on the $i$-th component $x_i:=\<x, e_i\>$, we have
 $P_t^0= \prod_{i=1}^\infty P_t^{0,i}$, so that this Harnack inequality leads to
 $$(P_t^0 f(x))^p\le P_t^0 f^p(y)  \exp\bigg[\ff{p}{2(p-1)} \sum_{i=1}^\infty  \ff{\ll_i|x_i-y_i|^2}{\e^{2\ll_i t}-1}\bigg],\ \ t>0, f\in \B_b^+(\H), x,y\in\H$$ for any $p>1$. Noting that $\mu_0$ is an invariant probability measure of $P_t^0$, by taking $p=2$ we obtain
 \beq\label{AB} (P_t^0 f(x))^2 \int_\H \exp\bigg[-\sum_{i=1}^\infty \ff{\ll_i(x_i-y_i)^2}{\e^{2\ll_i t}-1}\bigg]\mu_0(\d y) \le \mu_0(f^2),\ \ x\in\H, t>0, f\in L^2(\mu_0).\end{equation} Observing that
 $$\ff{\ll_i(x_i-y_i)^2}{\e^{2\ll_it}-1} +\ff{\ll_i y_i^2} 2 = \ff{\ll_i (\e^{2\ll_i t}+1)}{2(\e^{2\ll_it}-1)} \Big(y_i-\ff {2 x_i}{\e^{2\ll_i t}+1}\Big)^2 +
 \ff{\ll_ix_i^2}{\e^{2\ll_i t}+1},$$
 by \eqref{MU} we have
 \beg{align*} &\int_\H \exp\bigg[-\sum_{i=1}^\infty \ff{\ll_i(x_i-y_i)^2}{\e^{2\ll_i t}-1}\bigg]\mu_0(\d y)\\
 &= \prod_{i=1}^\infty \ff{\ss{\ll_i}}{\ss{2\pi}}
 \int_\R \exp\bigg[- \ff{\ll_i(x_i-y_i)^2}{\e^{2\ll_it}-1} -\ff{\ll_i y_i^2} 2 \bigg]\d y_i\\
 &=\exp\bigg[-\sum_{i=1}^\infty \ff{\ll_ix_i^2}{\e^{2\ll_i t}+1}\bigg]\bigg(\prod_{i=1}^\infty \ff{\e^{2\ll_i t}-1}{\e^{2\ll_it}+1}\bigg)^{\ff 1 2},
  \ \ t>0, x\in\H.\end{align*}
 So, \eqref{AB} reduces to
 \beq\label{AB1} P_t^0 f(x)\le \ss{\mu_0(f^2)}\, \GG_x(t),\ \ x\in\H, t>0, f\in L^2(\mu_0),\end{equation} where due to \eqref{EG},
\beq\label{AB1'}\beg{split}\GG_x(t)& := \exp\bigg[\ff 1 2\sum_{i=1}^\infty \ff{\ll_ix_i^2}{\e^{2\ll_i t}+1}\bigg]\bigg(\prod_{i=1}^\infty \ff{\e^{2\ll_i t}+1}{\e^{2\ll_it}-1}\bigg)^{\ff 1 4}\\
 &\le \exp\bigg[\ff 1 2\sum_{i=1}^\infty \ff{\ll_ix_i^2}{\e^{2\ll_i t}+1}\bigg] \bigg(\prod_{i=1}^\infty\Big(1+ \ff 1{\ll_i t}\Big)\bigg)^{\ff 1 4}<\infty, t>0, x\in\H.\end{split}\end{equation}
 Moreover, using the stronger condition $\sum_{i=1}^\infty \ll_i^{-\theta}<\infty$ for some $\theta\in (0,1)$ included in  {\bf (H3)}, and noting that $\log (1+ r)\le c r^\theta$ for some constant $c>0$ and all $r\ge 0$, we obtain
 \beq\label{AB3}\beg{split} &\Psi(t,x):=\int_0^t \log \GG_x(s) \d s = \ff 1 4 \sum_{i=1}^\infty\int_0^t \Big\{\ff{2\ll_i x_i^2}{\e^{2\ll_i s}+1} +\log \Big(1+\ff 1 {\ll_i s}\Big)\d s\Big\}\d s  \\
 &\le \ff 1 4 \sum_{i=1}^\infty \bigg\{2\int_0^t \ll_ix_i^2\e^{-2\ll_i s}\d s +\ff c{\ll_i^\theta} \int_0^{t} r^{-\theta} \d r\bigg\}\\
 &\le \ff 1 4\sum_{i=1}^\infty x_i^2(1-\e^{-2\ll_i t})    + C t^{1-\theta}  <\infty,\ \ t>0, x\in\H\end{split}\end{equation} for some constant $C>0.$ For later use we deduce from this that
 \beq\label{AB3'} \limsup_{t\to 0} \sup_{y\to x} \Psi(t,y) \le \ff 1 2 \limsup_{t\to 0} \sup_{y\to x}\bigg\{\sum_{i=1}^\infty x_i^2(1-\e^{-2\ll_i t}) +|x-y|^2   + C t^{1-\theta}\bigg\}=0.\end{equation}
 Since \eqref{AB1'} implies $\GG_x(s)\in (1,\infty)$, for every $t>0$ we have
 $$\bb_x(t):= \int_0^t \ff{\d s}{\GG_x(s)}\in (0,t],$$  so that
 $$\nu_{t,x}(\d s):= \ff{1_{[0,t]}(s)}{\bb_x(t) \GG_x(s)}\d s$$  is a probability measure on $[0,t].$
 Noting that $\ff{\bb_x(t)}t \int_0^t \GG_x(s)\nu_{t,x}(\d s)=1$ and $\log\big(\ff{\bb_x(t)}t \GG_x(s)\big)\le \log \GG_x(s)$, the Young inequality \eqref{Yang} yields
 \beg{align*} &\int_0^t |Z(X_s^x)|^2\d s = \ff{2t}\ll   \int_0^t \Big(\ff \ll 2 |Z(X_s^x)|^2\Big)\Big(\ff{\bb_x(t)}t\GG_x(s)\Big)\nu_{t,x}(\d s)\\
 &\le \ff{2 t}\ll \log\nu_{t,x}\big(\e^{\ff\ll 2 |Z(X_\cdot^x)|^2}\big) + \ff {2t} \ll \int_0^t \Big\{\ff{\bb_x(t)}t \GG_x(s)\log\Big(\ff{\bb_x(t)}t\GG_x(s)\Big)\Big\}\nu_{t,x}(\d s)\\
 &\le   \ff{2 t}\ll \log\nu_{t,x}\big(\e^{\ff\ll 2 |Z(X_\cdot^x)|^2}\big) + \ff 2 \ll \Psi(t,x),\ \ t\ge 0, x\in\H.\end{align*}
Combining this with \eqref{AB1} for $f= \e^{\ff\ll 2 |Z|^2}$,   \eqref{AB3} and $\mu_0(\e^{\ll |Z|^2})<\infty$, we arrive at
\beq\label{AB4} \beg{split} &\E\exp\bigg[\gamma \int_0^t |Z(X_s^x)|^2\d s\bigg] \le \e^{\ff{2\gamma}\ll \Psi(t,x)} \E\bigg\{\int_0^t \e^{\ff\ll 2 |Z(X_s^x)|^2}\nu_{t,x}(\d s)\bigg\}^{\ff{2\gamma t}\ll}\\
&\le \e^{\ff{2\gamma}\ll \Psi(t,x)} \bigg\{ \int_0^t \big\{P_s^0\e^{\ff \ll 2 |Z|^2}(x)\big\}\nu_{t,x}(\d s) \bigg\}^{\ff{2\gamma t}\ll} \\
&\le \e^{\ff{2\gamma}\ll \Psi(t,x)} \bigg\{ \int_0^t \ss{\mu_0(\e^{\ll|Z|^2})}\GG_x(s)\nu_{t,x}(\d s) \bigg\}^{\ff{2\gamma t}\ll} \\
&= \e^{\ff{2\gamma}\ll \Psi(t,x)}\Big(\ff t{\bb_x(t)} \ss{\mu_0(\e^{\ll|Z|^2})}\Big)^{\ff{2\gamma t}\ll} =:\LL(t,x,\gamma)<\infty,\ \ x\in \H, \gamma>0, t\in \Big(0, \ff\ll{2\gam}\Big].   \end{split}\end{equation}
By taking $\gamma=\ff 1 4,$ we prove \eqref{AB0} for   $t_0 =2\ll$.

\subsection{Strong Feller and strictly positive density of $P_t$}

By the Harnack inequality \eqref{HN}, $P_t^0$ is strong Feller having strictly positive density with respect to  $\mu_0$ (see \cite[Proposition 3.1(1)]{WY11}). Then as in (b) and (c) in the proof of Lemma \ref{P1},  we may prove the same property for $P_t$ using  \eqref{GR} and \eqref{AB4}. To save space,   we only prove here the strong Feller property.

For any $t>0,$ by the semigroup group property of $P_t$, \eqref{GR}, and the strong Feller property of $P_t^0$,    we obtain
\beq\label{AB5}\beg{split}&\limsup_{y\to x} |P_t f(x)-P_t f(y)|= \limsup_{r\to 0} \limsup_{y\to x} |P_r(P_{t-r} f)(x)-P_r(P_{t-r} f)(y)| \\
&\le  \limsup_{r\to 0} \limsup_{y\to x}\Big\{|P_r^0(P_{t-r} f)(x)- P_r^0(P_{t-r} f)(x)| \\
&\qquad \qquad \qquad\qquad +  \big| \E\big[(P_{t-r} f)(X_r^x)(R_r^x-1)- (P_{t-r} f)(X_r^y)(R_r^r-1)\big]\big|\Big\}\\
&\le \|f\|_\infty \limsup_{r\to 0} \limsup_{y\to x}\E\big(|R_r^x-1|+|R_r^y-1|\big).\end{split}\end{equation}
Recalling that $R_r^y= R_{0,r}^y$, by   \eqref{RX}  we have
$$\E|R_r^y-1|^2 = \E (R_r^y)^2-1\le \Big(\E \e^{3\int_0^r |Z(X_s^y)|^2\d s}\Big)^{\ff 1 2}-1,\ \ y\in\R^d.$$ So, according to \eqref{AB5},  $P_t$  is strong Feller provided
\beq\label{AB6} \limsup_{r\to 0} \limsup_{y\to x} \E \exp\bigg[3\int_0^r |Z(X_s^y)|^2\d s\bigg] =1.\end{equation}
Recall that $\bb_x(t)= \int_0^t \ff{\d s}{\GG_x(s)}\,\d s.$ By Jensen's inequality and \eqref{AB3} we have
$$\log\ff{\bb_x(t)}t =- \log \bigg(\ff 1 t \int_0^t \ff{\d s}{\GG_x(s)}\bigg) \le -\ff 1 t \int_0^t \Big\{\log \ff 1 {\GG_x(s)}\Big\}\d s= \ff{\Psi(t,x)}t.$$ Combining this with \eqref{AB3'} and \eqref{AB4}, we obtain
\beg{align*} \lim_{r\to 0}\limsup_{y\to x} \LL (r, y, 3)
&\le \lim_{r\to 0}\limsup_{y\to x} \e^{\ff 6 \ll \Psi(r,y)} \Big(\e^{\ff 1 r \Psi(r,y)} \ss{\mu_0(\e^{\ll|Z|^2})}\Big)^{\ff{6r}\ll}\\
&= \lim_{r\to 0}\limsup_{y\to x} \e^{\ff {12} \ll \Psi(r,y)}=1.\end{align*}
Combining this with \eqref{AB4}, we prove \eqref{AB6}.

\subsection{Uniqueness of the Feller semigroup and invariant probability measure}

  To prove that $P_t$ is the unique Feller Markov semigroup associated to \eqref{E2}, we recall the pathwise uniqueness for  $\mu_0$-a.e. initial points. By \cite[Theorem 1]{DFRV},  there exists an $\mu_0$-null set $\H_0$ such that for any $x\notin \H_0$, the SPDE \eqref{E2} has at most one mild solution starting at $x$ up to life time. Combining this with the weak solution constructed in (a), we see that for any initial point $x\notin \H_0$, the SPDE \eqref{E2} has a unique mild solution  $X_t^x$  which is non-explosive with distribution $P_t(x,\d y)$. So, if there exists another Feller  transition probability  kernel  $\bar P_t(x,\d y)$ associated to \eqref{E2}, then $\bar P_t(x,\d y)= P_t(x,\d y)$ for $x\notin \H_0.$ Since $\H\setminus \H_0$ is dense in $\H$, by the Feller property
these transition probability kernels are weak continuous in $x$, so that $\bar P_t(x,\d y)= P_t(x,\d y)$ for all $x\in \H$.

Next, according to \cite[Proposition 3.1(3)]{WY11}, to show that $P_t$ has at most one invariant probability measure, it suffices to prove for instance the Harnack inequality
\beq\label{HHL}(P_t f)^6(x)\le (P_t f^6)(y) H_t(x,y),\ \ x,y\in \H, f\in \B_b(\H)\end{equation} for some $t>0$ and measurable function $H_t: \H^2\to (0,\infty).$ By \eqref{GR} and \eqref{HN}, we have
\beg{align*} (P_tf(x))^6&= \big\{\E[f(X_t^x)R_t^x]\big\}^6 \le \big\{ P_t^0 f^2(x) \E (R_t^x)^2\big\}^4\\
&\le \big\{(P_t^0 f^4)(y) \big\}^2  \E(R_t^x)^6=  \big\{\E f^4(X_t^y) \big\}^2  [\E(R_t^x)^6]\\
&\le \big\{\E[f^6(X_t^y)R_t^y]\big\} \cdot \big\{\E (R_t^y)^{-1}\big\}\E(R_t^x)^6 = \{P_t f^6(y)\} \cdot \big\{\E (R_t^y)^{-1}\big\}[\E(R_t^x)^6].\end{align*}
By \eqref{AB4} and the definition of $R_t^\cdot$, it is easy to see that when $t>0$ is small enough,  $\{\E (R_t^y)^{-1} \}[\E(R_t^x)^6]\le H_t(x,y)$
holds for some measurable function $H_t: \H^2\to (0,\infty).$ Therefore, \eqref{HHL} holds.

\subsection{ $P_t$-invariance of $\mu$ and estimates on the density}
Finally, we prove that $\mu$  in Theorem \ref{T4.1} is an invariant probability measure of $P_t$. Let $\mu$ and $\bar \P_x$ be in Theorem \ref{T4.1}, according to the proof of Theorem \ref{T1.1}(2) we conclude that $\bar\P_\mu:=\int_\H \bar \P_x\mu(\d x)$ is the distribution of a weak solution to \eqref{E2} with initial distribution $\mu$. Since $\mu$ is absolutely continuous with respect to $\mu_0$, the uniqueness for $\mu_0$-a.e. initial points implies that the weak solution starting from $\mu$ is unique, so that $\mu(P_tf)=\mu(\bar P_t f)$ for $t\ge 0$ and $f\in \B_b(\H).$ Since $\mu$ is $\bar P_t$-invariant, it is $P_t$-invariant as well. Since Theorem \ref{T4.1} implies $\ss\rr \in H^{2,1}(\mu_0)$, \eqref{PD'} and \eqref{2.4'}, it remains to prove $\log\rr\in H^{2,1}(\mu_0)$.

By $\mu(\rr)=1$ and $\ss\rr\in H^{2,1}(\mu_0)$,  we have
$\log(\rr+\dd)\in H^{2,1}(\mu_0)$ for all $\dd>0$. Combining this with \eqref{2.4'} we conclude that $\log\rr\in H_\si^{2,1}(\mu_0)$ provided
$\mu_0(\rr>0)=1$ with $\mu_0(|\log\rr|^2)<\infty.$ It is well known that the Gaussian measure $\mu_0$ satisfies  the Poincar\'e inequality
$$\mu_0(f^2)\le \ff 1 {\ll_1}\mu_0(|\nn f|^2) + \mu_0(f)^2,\ \ f\in H_\si^{2,1}(\mu_0).$$
Then, as shown in the last step in the proof of Theorem \ref{T1.1}(1), $\mu_0(|\log\rr|^2)<\infty$ follows from \eqref{2.4'} if $\mu_0(\rr>0)=1$. Thus, we only need to prove $\mu_0(\rr>0)=1$.

Recalling that $R_t^x= R_{0,t}^x$ for $R_{s,t}^x$ defined in \eqref{RX}, by \eqref{GR} and \eqref{AB4} we may find out a constant $t_0>0$ and
some function $H\in C(\H\to (0,\infty))$ such that for any $f\in \B_b^+(\H)$,
\beg{align*}& \big(P_{t_0}^0 f(x)\big)^2 = \big(\E f(X_{t_0}^x)\big)^2 \le \big(\E [f^2(X_{t_0}^x)R_t^x]\big)  \E \big[(R_{t_0}^x)^{-2}\big]\\
&= (P_{t_0} f^2(x)) \E \big[(R_{t_0}^x)^{-2}\big]
\le H(x) P_{t_0} f^2(x),\ \ x\in\H.\end{align*}
Then for any measurable set $A\subset \H$ with $\mu_0(A)>0$, we have
\beq\label{HH}\mu(A) =\mu(P_{t_0}1_A^2) \ge \mu\bigg(\ff{(P_{t_0}^01_A)^2}{H}\bigg).\end{equation}
On the other hand,  by $\mu_0(P_{t_0}^0 1_A)=\mu_0(A)>0$, there exists $y\in \H$ such that $P_{t_0}1_A(y)>0$ so that \eqref{HN} implies
$$ P_{t_0}^0 1_A(x) \ge (P_{t_0}^0 1_A(y))^2 \e^{-\ff{C|x-y|^2}{t_0}}>0,\ \ x\in\H.$$ Combining this with \eqref{HH} and $\ff 1 {H}>0$.     Therefore, $\mu_0$ is absolutely continuous with respect to $\mu$ and hence, $\mu_0(\rr>0)=1.$

\section{Local Harnack inequality on incomplete manifolds}

Let $M$  be a $d$-dimensional differential manifold without boundary which is equipped with a (not necessarily complete) $C^2$-metric such that the curvature is well defined and continuous. Let $\DD$ and $\nn$ be the corresponding Laplace-Beltrami operator and the gradient operator respectively. Then for any $V\in C^2(M)$, the operator
$L:= \DD +\nn V$ generates a  unique diffusion process up to life time. Let $(X_t^x)_{t\in [0,\zeta(x)]}$ be the diffusion process starting at $x$ with life time $\zeta(x)$. Then the associated Dirichlet semigroup is given by
$$P_t f(x):= \E\big\{1_{\{t<\zeta(x)\}}f(X_t^x)\big\},\ \ x\in M, t\ge 0, f\in \B_b(M).$$
For any   $f\in \B_b^+(M):=\{f\in \B_b(M): f\ge 0\},$ define
$$E_{P_t}(f)= P_t (f\log f)- (P_tf)\log P_t f,\ \ t\ge 0.$$
Let $\rr$ be the Riemannian distance.  By the locally compact of the manifold  we may take $R\in C(M\to (0,\infty))$ such that
$$B_\rr(x, R(x)):=\big\{y\in M: \ \rr(x,y)\le R(x)\big\}$$ is compact for all $x\in M$. When the metric is complete this is true for all $R\in C(M\to (0,\infty)$.
We will use this function $R$ to establish the local Harnack inequality.

\beg{thm}\label{TA} There exists a function $H\in C(M\to (0,\infty))$ such that
\beq\label{AP0} |\nn P_t f(x) | \le \dd E_{P_t}(f)(x) + H(x) \Big(\dd + \ff 1 {\dd(t\land 1)}\Big),\ \ t>0, \dd\ge \ff{160}{R(x)}, f\in \B_b^+(M).\end{equation}
Consequently, for any $p>1$ there exists a function $F\in C(M\to (0,\infty))$ such that for any $t>0$ and  $f\in \B_b^+(M),$
\beq\label{AP0'}  (P_t f(x))^p \le (P_t f^p(y)) \exp\bigg[\ff{F(x)\rr(x,y)^2}{t\land 1}+ F(x)\bigg],
 \ \  x,y\in M\ \text{with}\ \rr(x,y)\le \ff 1 {F(x)}. \end{equation}
\end{thm}

\beg{proof} According to \cite{ATW06}, it is easy to prove  \eqref{AP0'}   from \eqref{AP0}. When the metric is complete, an estimate of type \eqref{AP0} for all $\dd>0$ has been proved in \cite{ATW09}. The only difference comes from the incompleteness of the metric for which we can not take $R(x)$ arbitrarily large as in \cite{ATW09}.   Below we   figure out the proof in the present case.

(1) To prove \eqref{AP0}, we fix      $f\in B_b^+(M)$. By using $\ff{f}{P_tf(x)}$ replace $f$, we may and do assume that $P_tf(x)=1$ at a fixed point $x$ so that $E_{P_t}(f)(x)=P_t(f\log f)(x).$

Now, let us check the proof of Theorem 1.1 in \cite{ATW09} (pages 3666-3667), where the part before (4.5) has nothing to do with the completeness; that is, with the compact set $D:= B_\rr(x, R(x))$, all estimates therein before (4.5) apply to the present setting. More precisely, letting
$$\tau(x)=\inf\{t\ge 0: X_t^x\notin D\},$$
we have ((4.1) in \cite{ATW09}) \beq\label{AP1} |\nn P_t f(x)|\le I_1+I_2,\end{equation}
where ((4.2) in \cite{ATW09})
\beq\label{AP2} I_1\le \dd \E\big\{1_{\{t<\tau(x)\}}(f\log f)(X_t^x)\big\} +\ff\dd \e + C(x)\Big(1+\ff 1 {\dd t}\Big),\ \ \dd>0, t>0\end{equation}
holds for function $C\in C(M\to (0,\infty))$ depending only on $d$ and curvature of the operator $L$; and moreover ((4.5) in \cite{ATW09}),
\beq\label{AP3} I_2\le \dd \E\big\{1_{\{\tau(x)\le t<\zeta(x)\}}(f\log f)(X_t^x)\big\} +\ff\dd \e + \dd \log \E \e^{\ff{9R(x)}{\dd\tau(x)}} +A(x),\ \ \dd>0, t>0\end{equation}
holds for $A(x):= \sup_{r>0} \big\{C(x)\ss r\log (\e+r) -r\big\}$, which is finite and continuous in $x$.  Now, due to the restriction of $R(x)$, we have to take large enough $\dd>0$ and can not replace $\dd$ by $\dd\land 1$ as in (4.5) of \cite{ATW09}. This will lead to less harp estimate but it is enough for our study in the present paper. More precisely, using $\dd$ to replace $\aa\land 1$ in the display after (4.5) of \cite{ATW09}, we have
$$\E \e^{\ff{9R(x)}{\dd\tau(x)}} \le 1 + 9 \int_0^\infty (9u+1)\e^{-u}\d u=:A'<\infty,\ \ \dd\ge \ff{160}{R(x)}.$$
Combining the with \eqref{AP1}-\eqref{AP3}, we prove \eqref{AP0} for some $H\in C(M\to (0,\infty)).$

(2) Since $H,R$ are strictly positive and continuous, and $B_\rr(x,R(x))$ is compact for every $x$,
$$ \bar H(x):= \sup_{B_\rr(x,R(x))}H \ \ \text{and}\ \  \hat R(x):= \inf_{B_\rr(x,R(x))}R $$ are strictly positive continuous functions in $x$. For any $p>1$, let
$$G(x)= \ff{p-1}{p\bar H(x)}\land\hat R(x),\ \ x\in M.$$ Then \eqref{AP0} implies
$$ |\nn P_t f(y)|\le \dd E_{P_t}(f)(y) +\bar H(x) \Big(\ff 1 {\dd(1\land t)}+\dd\Big),\ \ y\in B_\rr (x, G(x)), \dd\ge \ff{160}{\hat R(x)}  $$ for
$f\in \B_b^+(M).$
So, letting $\gam: [0,1]\to M$ be the minimal geodesic from $x$ to $y$ with $|\dot \gam_s|= \rr(x,y)$ for $s\in [0,1],$   letting $\bb(s)= 1+s(p-1)$, and applying  the above inequality with $\dd:=\ff{p-1}{p\rr(x,y)}\ge \ff{160}{\hat R(x)}$, we obtain
\beg{align*} &\ff{\d}{\d s} \big\{\log P_t f^{\bb(s)}\big\}^{\ff p {\bb(s)}} = \ff{p(p-1) E_{P_t}(f^{\bb(s)})}{\bb(s)^2 P_t f^{\bb(s)}} + \ff{p \<\nn P_t f^{\bb(s)},\dot\gam_s\>}{\bb(s) P_t f^{\bb(s)}} (\gam_s)\\
&\ge \ff{p\rr(x,y)}{\bb(s)P_t f^{\bb(s)}(\gam_s)}\bigg\{ \ff{p-1}{p\rr(x,y)} E_{P_t}(f^{\bb(s)}) - |\nn P_t f^{\bb(s)}|\bigg\}(\gam_s)\\
&\ge -\ff{p\rr(x,y)}{\bb(s) P_t f^{\bb(s)}(\gam_s)} \bigg\{\bar H(x) \big(P_t f^{\bb(s)}(\gam_s)\big) \Big(\ff{p\rr(x,y)}{(p-1)(t\land 1)} +\ff{p-1}{p\rr(x,y)}\Big)\bigg\}\\
&\ge -  \bar H(x) \Big(\ff{p^2\rr(x,y)^2}{(p-1)(t\land 1)}+1\Big),\ \ s\in [0,1],\ \rr(x,y)\le G(x).\end{align*} Integrating over $[0,1]$ with respect to $\d s$, we prove \eqref{AP0'} for $F:= \ff{p^2\bar H}{p-1}\lor \ff 1 {G}.$

\end{proof}

\paragraph{Acknowledgements.} The author would like to   thank   Vladimir Bogachev, Chenggui Yuan and the referees for helpful comments and corrections.

\beg{thebibliography}{99}

\bibitem{ABR} S. Albeverio, V. Bogachev, M. R\"ockner, \emph{On uniqueness of invariant measures for finite-dimensional diffusions,} Comm. Pure Appl. Math. 52(1999), 325--362.
 \bibitem{ATW06} M. Arnaudon, A. Thalmaier, F.-Y. Wang, \emph{ Harnack inequality and heat kernel estimates
  on manifolds with curvature unbounded below,} Bull. Sci. Math. 130(2006), 223--233.
\bibitem{ATW09} M. Arnaudon, A. Thalmaier, F.-Y. Wang, \emph{Gradient estimates and Harnack inequalities on non-compact Riemannian manifolds,} Stoch. Proc. Appl. 119(2009), 3653--3670.

\bibitem{ATW13}  M. Arnaudon, A. Thalmaier, F.-Y. Wang, \emph{Equivalent log-Harnack and gradient for point-wise curvature lower bound,}  Bull. Math. Sci. 138(2014), 643--655.

\bibitem{BE}   D. Bakry, D. and M. Emery, \emph{Hypercontractivit\'e de
semi-groupes de diffusion}, C. R. Acad. Sci. Paris. S\'er. I Math.
299(1984), 775--778.


\bibitem{BC} R. F. Bass, Z. Q. Chen, \emph{ Brownian motion with singular drift,}  Ann. Probab. 31(2003), 791--817.

\bibitem{BKR0} V. I. Bogachev, N. V. Krylov,   M. R\"ockner, \emph{On regularity of transition probabilities
and invariant measures of singular diffusions under minimal conditions,} Comm. Part.
Diff. Equat. 26 (2001), 2037--2080.
\bibitem{BKR}   V.I.~Bogachev,  N.V.~Krylov,  M.~R\"ockner, \emph{
 Elliptic and parabolic equations for measures,} Russ. Math. Surv. 64(2009), 973--1078.

\bibitem{[5]}  V.I.~Bogachev, M.~R\"ockner,
 \emph{Regularity of invariant
measures on finite and infinite dimensional spaces and applications,}
  J. Funct. Anal.  133(1995), 168--223.

\bibitem{[2]}  V.I.~Bogachev,  M.~R\"ockner,
\emph{  A generalization
of Hasminskii's theorem on existence of invariant measures for locally
integrable drifts,}  Theo. Probab. Appl. 45(2002), 363--378.

\bibitem{BR99}
  V.I.~Bogachev,  M.~R\"ockner,
\emph{ Elliptic equations for measures on infinite dimensional spaces
and applications,}   Probab. Theo. Relat. Fields 120(2001), 445--496.

\bibitem{BRW}	V. I. Bogachev, M. R\"ockner, F.-Y. Wang, \emph{Elliptic equations for invariant measures on finite and infinite dimensional manifolds,} J. Math. Pure Appl. 80(2001), 177--221.

\bibitem{CGW}
 P. Cattiaux,  A. Guillin,  L. Wu, \emph{A note on Talagrand's transportation inequality and logarithmic Sobolev
 inequality,} Probab. Theo. Relat. Fields 148(2010), 285--304.

\bibitem{CDR}  P. E. Chaudru de Raynal, \emph{Strong existence and uniqueness for stochastic differential equation with H\"orlder drift and degenerate noise,} to appear in  {\it Ann. Inst. Henri Poincar\'e Probab. Stat.}
http://arxiv.org/abs/1205.6688.


 \bibitem{DFRV}   G. Da Prato, F. Flandoli, M. R\"ockner,  A. Yu. Veretennikov, \emph{Strong uniqueness for stochastic evolution equations
with unbounded measurable drift term,} J. Theo. Probab. 28(2015), 1571--1600.

 \bibitem{DR} G. Da Prato, M. R\"ockner, \emph{Singular dissipative stochastic equations in Hilbert spaces,}
Probab. Theo. Relat. Fields 124(2002), 261--303.

\bibitem{DZ}  G. Da Prato, J. Zabczyk,  \emph{Stochastic Equations in Infinite Dimensions,}  Cambridge University Press, Cambridge, 1992.

\bibitem{DS}   E. B. Davies, B. Simon, \emph{Ultracontractivity and heat kernel
for Schr\"odinger operators and Dirichlet Laplacians,} J. Funct. Anal. 59(1984), 335--395.

\bibitem{PW} E. Priola, F.-Y. Wang, \emph{Gradient estimates for diffusion semigroups with singular coefficients, } J. Funct. Anal. 236(2006), 244--264.

\bibitem{GW} F. Gong, F.-Y. Wang,  \emph{Functional inequalities for uniformly integrable semigroups and application to essential spectrum,} Forum Math. 14(2002), 293--313.

\bibitem{Gross}  L. Gross,  \emph{Logarithmic Sobolev inequalities,} Amer. J. Math.
  97(1976), 1061--1083.

\bibitem{GM} I. Gy\"ongy, T. Martinez, \emph{On stochastic differential equations with locally unbounded drift,}
Czechoslovak Math. J. 51(2001), 763--783.

\bibitem{H80}  R.Z.~Hasminskii, \emph{Stochastic Stability of Differential Equations,} Sijthoff and
noordhoff, 1980.

\bibitem{IW} N. Ikeda, S. Watanabe, \emph{Stochastic
Differential Equations and Diffusion Processes (Second Edition),}
  North-Holland, 1989.

\bibitem{KR}   N. V. Krylov, M. R\"ockner, \emph{Strong solutions of stochastic equations with singular time dependent drift,}  Probab. Theo. Relat. Fields   131(2005), 154--196.

\bibitem{RW} 	M. R\"ockner, F.-Y. Wang, \emph{Supercontractivity and ultracontractivity for (non-symmetric) diffusion semigroups on manifolds,}
 Forum Math.  15(2003), 893--921.

 \bibitem{RW03}  M. R\"ockner, F.-Y. Wang, \emph{Harnack and functional inequalities for generalized Mehler semigroups,} J. Funct. Anal. 203(2003), 197--234.
  \bibitem{RW03b}  M. R\"ockner, F.-Y. Wang, \emph{ Supercontractivity and ultracontractivity for (non-symmetric) diffusion semigroups on manifolds,} Forum Math. 15(2003), 893--921.

\bibitem{ST} W. Stannat,  \emph{(Nonsymmetric) Dirichlet operators on $L^1$: Existence, uniqueness and associated Markov processes,}
Ann. Sc. Norm. Super. Pisa Cl. Sci. (4) 28(1999) 99--140.

\bibitem{SV} D. Stroock, S.R.S. Varadhan, \emph{Multidimensional Diffusion processes,} Springer, 1979.

\bibitem{TT}  M. Takeda, G. Trutnau, \emph{ Conservativeness of non-symmetric diffusion processes generated by perturbed divergence forms,}  Forum Math. 24(2012),   419--444.


\bibitem{W97}  F.-Y. Wang, \emph{Logarithmic Sobolev inequalities on noncompact Riemannian manifolds,} Probab. Theo. Relat. Fields 109(1997), 417--424.

\bibitem{W00a} F.Y.~Wang, \emph{Functional inequalities for empty essential spectrum,}
 J. Funct. Anal.  170 (2000), 219--245.

\bibitem{W00b} F.-Y. Wang, \emph{ Functional inequalities, semigroup properties and spectrum estimates,}  Inf. Dimens. Anal. Quant. Probab. Relat. Top. 3(2000), 263--295.

\bibitem{W05} F.-Y. Wang, \emph{Functional Inequalities, Markov Processes and Spectral Theory,} Science Press,   2005.

\bibitem{W01}  F.-Y. Wang, \emph{Logarithmic Sobolev inequalities: conditions and counterexamples, } J. Operat. Theo. 46(2001), 183--197.

\bibitem{Wbook}  F.-Y. Wang, \emph{Analysis for Diffusion Processes on Riemannian Manifolds,} World Scientific, 2014.

\bibitem{Wbook2} F.-Y. Wang, \emph{Harnack Inequalities for Stochastic Partial Differential Equations,} Springer, 2013.

\bibitem{W15}    F.-Y. Wang, \emph{Gradient estimates and applications for SDEs in Hilbert space with multiplicative
noise and Dini drift,}  J. Diff. Equat. 260(2016), 2792--2829.

 \bibitem{WY11} 	F.-Y. Wang, C. Yuan, \emph{Harnack inequalities for functional SDEs with multiplicative noise and applications,} Stoch. Proc. Appl. 121(2011), 2692--2710.

\bibitem{WZa} F.-Y. Wang, X. Zhang, \emph{Degenerate SDEs in Hilbert spaces with rough drifts,}   Infin. Dimens. Anal. Quantum Probab. Relat. Top.  18(2015), 1550026, 25 pp.

\bibitem{WZb}  F.-Y. Wang, X. Zhang, \emph{Degenerate SDE with H\"older-Dini drift and non-Lipschitz noise coefficient,} to appear in SIMA J. Math. Anal.

\bibitem{Wu} L.-M.  Wu, \emph{ Moderate deviations of dependent random variables related to CLT,}   Ann. Probab. 23(1995), 420--445.

\bibitem{Zhang} X. Zhang, \emph{ Strong solutions of SDES with singular drift and Sobolev diffusion coefficients,} \emph{Stoch. Proc. Appl. } 115(2005), 1805--1818.
\bibitem{Zhang2} X. Zhang, \emph{Stochastic homeomorphism flows of SDEs with singular drifts and Sobolev diffusion coefficients,} Electr. J. Probab.
 16(2011), 1096--1116.
\end{thebibliography}
\end{document}